\documentclass[review]{elsarticle}
\usepackage[fleqn]{amsmath}
\usepackage{graphicx}
\usepackage{subfigure}
\usepackage{lineno,hyperref}
\modulolinenumbers[5]
\usepackage{CJK}
\usepackage{indentfirst}
\usepackage{amsmath}

\journal{\textbf{XXX}}

\begin{document}

\begin{frontmatter}


\title{On the study of stochastic fractional-order differential equation systems}

\author[mymainaddress]{Guang-an Zou\corref{mycorrespondingauthor}}
\cortext[mycorrespondingauthor]{Corresponding author}
\ead{zouguangan00@163.com}

\author[mymainaddress]{Bo Wang}

%\author{Guang-an Zou\fnref{myfootnote}}
%\address{School of Mathematics and Statistics, Henan University, Kaifeng 475004, P. R. China}


\address[mymainaddress]{School of Mathematics and Statistics, Henan University, Kaifeng 475004, China}
%\address[mysecondaryaddress]{Institute of Applied Mathematics, Henan University, 475004 Kaifeng, China}

\begin{abstract}

In this article, the existence and uniqueness about the solution for a class of stochastic fractional-order differential equation systems are investigated, where the fractional derivative is described in Caputo sense. The fractional calculus, stochastic analysis techniques and the standard Picard's iteration are used to obtain the required results, the nonlinear term is satisfied with some non-Lipschitz conditions (where the classical Lipschitz conditions are special cases). The stochastic fractional-order Newton-Leipnik and Lorenz systems are provided to illustrate the obtained theory, and numerical simulation results are also given by the modified Adams predictor-corrector scheme.

\end{abstract}

\begin{keyword}
Fractional-order differential equations, stochastic Newton-Leipnik system, stochastic Lorenz system, numerical simulations.
%\texttt{elsarticle.cls}\sep \LaTeX\sep Elsevier \sep template
%\MSC[2010] 00-01\sep  99-00
\end{keyword}

\end{frontmatter}

%\linenumbers

\section{Introduction}

Fractional calculus and fractional-order differential equations have been widely applied in many fields of science and engineering, such as physics [1-2], chemical [3-5], mechanics [6-8], biological [9-10], medical [11-13], economics and finance [14-15], materials [16-17], control theory [18-19], etc. Actually, the concepts of fractional derivatives are not only generalization of the ordinary derivatives, but also it has been found that they can efficiently and properly describe the behavior of many physical systems (real-life phenomena) more accurately than integer order derivatives. The advantages of fractional derivatives are that they have a greater degree of flexibility in the model and provide an excellent instrument for the description of the ``memory'' and ``hereditary'' properties of various practical processes and dynamical systems, which take into account the fact that the future state depend not just upon the present states but also upon all of the history of its previous states [20]. In particular, the chaotic dynamical systems described by fractional differential equations have gained more considerable and attentions. Some examples of the chaotic systems include the fractional-order Lorenz system, Chen system, Chua system, L\"{u} system, Liu system, financial system and Newton-Leipnik system have been studied in the literature [15,21-24]. Furthermore, theories on the existence, uniqueness and stability of solutions of initial-value problems for fractional-order differential equations have been established [25-28].

As the stochastic disturbances are unavoidable, in recent years, stochastic differential equations have become more and more important and interesting to researchers due to their successful and potential applications in various fields [29-31], and the basic theories and results of stochastic differential equations can be found in [29]. Studying of stochastic dynamical systems has been carried out by various researchers. For examples, the random attractor and stochastic bifurcation behavior of the stochastic Lorenz system have studied in [32-33]. Anishchenko et al. have investigated the phenomenon of stochastic resonance for chaotic systems perturbed by white noise and a harmonic force [34]. In Ref [35], it has been shown that the chaotic transitions in stochastic dynamical systems. In addition, the analytic and numerical methods for stochastic dynamical systems also have been constructed by some authors [36-37]. However, there are relatively few studies of stochastic differential equations with fractional derivatives. The existence and uniqueness for a class of stochastic fractional-order differential equations were obtained in [38-39]. The approximate controllability of fractional stochastic dynamic systems has been proved in [40]. But there has been little mention of the dynamic systems described by stochastic fractional-order differential equation systems. Motivated by the above facts, in this paper, we establish the existence and uniqueness of the solution for a class of stochastic fractional-order differential equation systems, and give the numerical results of stochastic fractional-order dynamic systems.

The paper is organized as follows: In Section 2, we present some essential facts in fractional calculus and stochastic analysis that will be used to obtain our main results. In Section 3, the existence and uniqueness of the solution for stochastic fractional-order differential equation systems are proved by Picard's iteration. In Section 4, the stochastic fractional-order Newton-Leipnik and Lorenz systems are presented to illustrate the obtained theory, and numerical simulation results are also given. Finally, the conclusions are drawn in Section 5.

\section{Notations and preliminaries}

In this section, we give some basic definitions, notations and lemmas which will be used throughout the paper, in order to establish our main results.

First of all, we define the infinite-dimensional space $\ell^{2}=\{\mathbf{x}=(x_{i})_{i\in Z},x_{i}\in R:\sum\limits_{i\in Z}x_{i}^2 <+\infty\}$ with the inner product and norm:
\begin{align*}
(\mathbf{x},\mathbf{y})=\sum_{i\in Z} x_{i}y_{i}, \|\mathbf{x}\|^{2}=\sum_{i\in Z} x_{i}^{2}, \forall~\mathbf{x}=(x_{i})_{i\in Z}, \mathbf{y}=(y_{i})_{i\in Z}\in \ell^{2}.\tag{2.1}
\end{align*}
which is a Hilbert space.

Let $(\Omega,\mathcal{F},\mathbf{P})$ be a complete probability space, for a separable Hilbert space $H$ with inner product $(\cdot,\cdot)$ and norm $\|\cdot\|$. Then $L_{2}(\Omega,H)$ is Hilbert space of $H$-valued random variables with the inner product $\mathbf{E}(\cdot,\cdot)$ and the norm $(\mathbf{E}\|\cdot\|^{2})^{1/2}$, in which $\mathbf{E}$ denotes the expectation.

For $v\in L_{2}(\Omega,H)$, there holds the following It\^{o} isometry property:
\begin{align*}
\mathbf{E}\|\int_{0}^{t}v(s)dW(s)\|^{2}=\int_{0}^{t}\mathbf{E}\|v(s)\|^{2}ds. \tag{2.2}
\end{align*}
where ${W(t)}_{t\geq 0}$ is a Wiener process.

Secondly, let us introduce three common notation for the fractional-order differential operator: the Riemann-Liouville, the Caputo-type, and the Gr\"{u}nwald-Letnikov fractional derivative.
For more details see [1,28].

\textbf{Definition 1.} The Riemann-Liouville fractional derivative of $f$ is defined as
\begin{align*}
^{R}D_{t}^{\alpha}f(t)=\frac{1}{\Gamma(n-\alpha)}\frac{d^{n}}{dt^{n}}\int_{0}^{t}\frac{f(s)}{(t-s)^{\alpha+1-n}}ds, t>0,n-1<\alpha<n,  \tag{2.3}
\end{align*}
where $\Gamma(\cdot)$ stands for the gamma function $\Gamma(x)=\int_{0}^{\infty}t^{x-1}e^{-t}dt$, and $n=[\alpha]+1$ with $[\alpha]$ denotes the integer part of $\alpha$.

\textbf{Definition 2.} The Caputo-type derivative of order $\alpha$ for a function $f$ can be written as
\begin{align*}
^{C}D_{t}^{\alpha}f(t)=\frac{1}{\Gamma(n-\alpha)}\int_{0}^{t}\frac{f^{n}(s)}{(t-s)^{\alpha+1-n}}ds, t>0, n-1<\alpha<n.\tag{2.4}
\end{align*}

\textbf{Definition 3.} The Gr\"{u}nwald-Letnikov fractional derivative of $f$ is given by
\begin{align*}
^{G}D_{t}^{\alpha}f(t)=\mathop{\mathrm{lim}}_{h\rightarrow 0}\frac{1}{h^{\alpha}}\mathop{\Sigma}_{j=0}^{[\frac{t-\alpha}{h}]}(-1)^{j}\tbinom{\alpha}{j}f(t-jh), \alpha>0. \tag{2.5}
\end{align*}

\textbf{Remark 2.1.}(1) The relationship between the Riemann-Liouville derivative and the Caputo-type derivative can be written as
\begin{align*}
^{C}D_{t}^{\alpha}f(t)=^{R}D_{t}^{\alpha}[f(t)-\mathop{\Sigma}_{k=0}^{n-1}\frac{t^{k}}{k!}f^{(k)}(0)].\tag{2.6}
\end{align*}

(2) The Riemann-Liouville derivatives is often approximated by using the Gr\"{u}nwald-Letnikov definition based on finite differences [41].

(3) The Caputo-type derivative of a constant is equal to zero.

In this study, we consider the Caputo-type fractional derivative of order $\alpha$ for a vector-valued function $\mathbf{y}(t)$, and the initial value problem of stochastic fractional-order differential equation is given as following:
\begin{align*}
\begin{cases}
^{C}D_{t}^{\alpha}\mathbf{y}(t)=f(t,\mathbf{y}(t))+\sigma(t,\mathbf{y}(t))\dot{W}(t),0\leq t\leq T,\\
\mathbf{y}^{(k)}(0)=\mathbf{y}_{0}^{(k)}, k=0,1,2,\ldots,m-1,
\end{cases}  \tag{2.7}
\end{align*}
where the functions $f(t,\mathbf{y}(t))$ and $\sigma(t,\mathbf{y}(t))$: $[0,T]\times R^{d}\rightarrow R^{d}$ are vector field, and the dimension $d\geq1$. The term $\dot{W}(t)=\frac{dW}{dt}$ describes a state dependent random noise, ${W(t)}_{t\geq 0}$ is a a standard scalar Brownian motion or Wiener process defined on a given filtered probability space $(\Omega,\mathcal{F},\{\mathcal{F}_{t}\}_{t\geq 0},\mathbf{P})$ with a normal filtration $\{\mathcal{F}_{t}\}_{t\geq 0}$, which is an increasing and continuous family of $\sigma$-algebras of $\mathcal{F}$, contains all of $\mathbf{P}$-null sets, and $W(t)$ is $\mathcal{F}_{t}$-measurable for each $t\geq 0$.

Here, let us recall the definitions of fractional calculus [28], the fractional integral operator of order $\alpha$ is given as following
\begin{align*}
I^{\alpha}\mathbf{g}(t)=\frac{1}{\Gamma(\alpha)}\int_{0}^{t}(t-s)^{\alpha-1}\mathbf{g}(s)ds, t>0. \tag{2.8}
\end{align*}

Applying the integral operator (2.8) to the both sides of initial value problem (2.7), we can obtain the Volterra integral equation
\begin{align*}
\mathbf{y}(t)=&\mathop{\Sigma}_{k=0}^{[\alpha]-1}\frac{t^{k}}{k!}\mathbf{y}^{(k)}(0)+\frac{1}{\Gamma(\alpha)}\int_{0}^{t}(t-s)^{\alpha-1}f(s,\mathbf{y}(s))ds\\
&+\frac{1}{\Gamma(\alpha)}\int_{0}^{t}(t-s)^{\alpha-1}\sigma(s,\mathbf{y}(s))dW(s) ,\tag{2.9}
\end{align*}
where $n-1<\alpha<n$ and $t\geq0$. Conversely, substituting the fractional derivative (2.4) into the equation (2.9), its leads to initial value problem (2.7). Consequently, we have the following Lemma.

\textbf{Lemma 2.1.} Every solution of the Volterra integral equation (2.9) is also a solution of the original initial value problem (2.7), and vice versa.

Therefore, we may focus our attention on equation (2.9), in other words, this allows us only discuss the properties of the solution of the equation (2.9) instead of the initial value problem (2.7). However, when the fractional order $\alpha\in(0,1)$, the equation (2.9) is singular but regular when $\alpha\geq1$. In the situation of singular, the Volterra equation (2.9) can be written as
\begin{align*}
\mathbf{y}(t)=&\mathbf{y}_{0}+\frac{1}{\Gamma(\alpha)}\int_{0}^{t}(t-s)^{\alpha-1}f(s,\mathbf{y}(s))ds\\
&+\frac{1}{\Gamma(\alpha)}\int_{0}^{t}(t-s)^{\alpha-1}\sigma(s,\mathbf{y}(s))dW(s) , \tag{2.10}
\end{align*}
in which $\alpha\in(0,1)$ and $t\geq0$.

\section{Existence and uniqueness result}

In this section, we will study the existence and uniqueness of solution to equation (2.10). Throughout the paper the following non-Lipschitz conditions are assumed and imposed:

(A1) The functions $f$ and $\sigma$ are measurable and continuous in $H$ for each fixed $t\in [0,T]$ and there exists a bounded function $L$: $[0,T]\times [0,+\infty]\rightarrow [0,+\infty]$, $(t,u)\rightarrow L(t,u)$ such that
\begin{align*}
 \mathbf{E}(\|f(t,\mathbf{x})\|^{2})+\mathbf{E}(\|\sigma(t,\mathbf{x})\|^{2})\leq L(t,\mathbf{E}(\|\mathbf{x}\|^{2})),  \tag{3.1}
\end{align*}
for all $t\in R$ and all $\mathbf{x}\in L_{2}(\Omega,H)$.

(A2) There exists a bounded function $K$: $[0,T]\times [0,+\infty]\rightarrow [0,+\infty]$ such that
\begin{align*}
\mathbf{E}(\|f(t,\mathbf{x})-f(t,\mathbf{y})\|^{2})+\mathbf{E}(\|\sigma(t,\mathbf{x})-\sigma(t,\mathbf{y})\|^{2})\\
\leq K(t,\mathbf{E}(\|\mathbf{x}-\mathbf{y}\|^{2})), \tag{3.2}
\end{align*}
for all $t\in R$ and all $\mathbf{x},\mathbf{y}\in L_{2}(\Omega,H)$.

\textbf{Remark 3.1.} (1) If the function $K(t,u)=Mu,u\geq 0$, and $M>0$ is a constant, then the condition (A2) implies global Lipschitz condition.

(2) If $K(t,u)$ is concave for each fixed $t\geq0$, and
\begin{align*}
\|f(t,\mathbf{x})-f(t,\mathbf{y})\|^{2}+\|\sigma(t,\mathbf{x})-\sigma(t,\mathbf{y})\|^{2}
\leq K(t,\|\mathbf{x}-\mathbf{y}\|^{2}), \tag{3.3}
\end{align*}
for all $\mathbf{x},\mathbf{y}\in H$ and $t\geq0$, the condition (A2) is still satisfied by Jensen's inequality.

(3) If $K(t,u)=m(t)\rho(u)$ with $u\geq0$ and $m(t)\geq0$ is locally integrable, and $\rho$ is continuous, monotone non-decreasing and concave function with $\rho(0)=0$, $\rho(u)>0$ for $u>0$ and $\int 1/\rho(u)du=\infty$, then the $\rho(u)$ satisfies the condition (A2).

Obviously, the classical Lipschitz conditions are only special cases in the above descriptions.

\textbf{Lemma 3.1.} If the function $L(t,u)$ is locally integrable in $t$ for each fixed $u\in[0,+\infty)$ and is continuous non-decreasing in $u$ for each fixed $t\in [0,T]$, for all $\lambda>0$, $u_{0}\geq0$, then the integral equation
\begin{align*}
u(t)=u_{0}+\lambda\int_{0}^{t}L(s,u(s))ds, \tag{3.4}
\end{align*}
has a global solution on $[0,T]$.

\textbf{Lemma 3.2.} The function $K(t,u)$ is locally integrable in $t$ for each fixed $u\in[0,+\infty)$ and is continuous non-decreasing in $u$ for each fixed $t\in [0,T]$, for $K(t,0)=0$ and $\gamma>0$, if a non-negative continuous function $z(t)$ satisfies
\begin{align*}
\begin{cases}
z(t)\leq \gamma \int_{0}^{t}K(s,z(s))ds,t\in R,\\
z(0)=0,
\end{cases} \tag{3.5}
\end{align*}
then $z(t)=0$ for all $t\in [0,T]$.

In order to consider the existence and uniqueness of the solution of equation (2.10), we attempt to use the following approximate technique, known as Picard's iteration. The sequence of stochastic process $\{\mathbf{y}_{n}\}_{n\geq0}$ is constructed as follows:
\begin{align*}
\begin{cases}
\mathbf{y}_{0}(t)=\mathbf{y}_{0},\\
\mathbf{y}_{n+1}(t)=\mathbf{y}_{0}+G_{1}(\mathbf{y}_{n})(t)+G_{2}(\mathbf{y}_{n})(t),n\geq1,
\end{cases}\tag{3.6}
\end{align*}
in which
\begin{align*}
\begin{cases}
G_{1}(\mathbf{y}_{n})(t)=\frac{1}{\Gamma(\alpha)}\int_{0}^{t}(t-s)^{\alpha-1}f(s,\mathbf{y}_{n}(s))ds,\\
G_{2}(\mathbf{y}_{n})(t)=\frac{1}{\Gamma(\alpha)}\int_{0}^{t}(t-s)^{\alpha-1}\sigma(s,\mathbf{y}_{n}(s))dW(s).
\end{cases}\tag{3.7}
\end{align*}

Take into account the proof of the existence and uniqueness result, we need the following two lemmas.

\textbf{Lemma 3.3.} The sequence of stochastic process $\{\mathbf{y}_{n}\}_{n\geq0}$ in (3.6) is bounded in $L_{2}(\Omega,H)$, i.e., $\mathrm{sup}_{n\geq0}\|\mathbf{y}_{n}\|_{L_{2}(\Omega,H)}\leq C$, where $C$ is a constant.

\textbf{Proof.} Via the inequality
\begin{align*}
(a+b+c)^{n}\leq 3^{n-1}(a^{n}+b^{n}+c^{n}), a,b,c\geq 0,n\geq1,\tag{3.8}
\end{align*}
we have
\begin{align*}
\mathbf{E}\|\mathbf{y}_{n+1}(t)\|^{2}\leq 3\mathbf{E}\|\mathbf{y}_{0}\|^{2}+3\mathbf{E}\|G_{1}(\mathbf{y}_{n})(t)\|^{2}+3\mathbf{E}\|G_{2}(\mathbf{y}_{n})(t)\|^{2}.\tag{3.9}
\end{align*}

Using the H\"{o}lder's inequality and the assumptions (3.1) and $\alpha>1/2$ for the right hand side of the above inequality, we can obtain
\begin{align*}
\mathbf{E}\|G_{1}(\mathbf{y}_{n})(t)\|^{2}&\leq\frac{1}{\Gamma^{2}(\alpha)}\mathbf{E}\|\int_{0}^{t}(t-s)^{\alpha-1}f(s,\mathbf{y}_{n}(s))ds\|^{2}\\
&\leq\frac{t^{2\alpha-1}}{\Gamma^{2}(\alpha)(2\alpha-1)}\int_{0}^{t}\mathbf{E}(\|f(s,\mathbf{y}_{n}(s))\|^{2})ds\\
&\leq k_{1}\int_{0}^{t}L(s,\|\mathbf{y}_{n}(s)\|^{2}_{L_{2}(\Omega,H)})ds, \tag{3.10}
\end{align*}
where $k_{1}=\frac{T^{2\alpha-1}}{\Gamma^{2}(\alpha)(2\alpha-1)}$.

Applying the It\^{o} isometry property (2.2), the H\"{o}lder's inequality and the assumptions (3.1) and $\alpha>1/2$ to the right hand side of the inequality (3.9), we have
\begin{align*}
\mathbf{E}\|G_{2}(\mathbf{y}_{n})(t)\|^{2}&\leq\frac{1}{\Gamma^{2}(\alpha)}\int_{0}^{t}\mathbf{E}\|(t-s)^{\alpha-1}\sigma(s,\mathbf{y}(s))\|^{2}ds\\
&\leq\frac{t^{2\alpha-1}}{\Gamma^{2}(\alpha)(2\alpha-1)}\int_{0}^{t}\mathbf{E}(\|\sigma(s,\mathbf{y}_{n}(s))\|^{2})ds\\
&\leq k_{1}\int_{0}^{t}L(s,\|\mathbf{y}_{n}(s)\|^{2}_{L_{2}(\Omega,H)})ds. \tag{3.11}
\end{align*}

Therefore, using the above relations (3.10) and (3.11) into the estimate (3.9), we have
\begin{align*}
\|\mathbf{y}_{n+1}(t)\|^{2}_{L_{2}(\Omega,H)}\leq C_{1}+C_{2}\int_{0}^{t}L(s,\|\mathbf{y}_{n}(s)\|^{2}_{L_{2}(\Omega,H)})ds,\tag{3.12}
\end{align*}
in which $C_{1}=3\mathbf{E}\|\mathbf{y}_{0}\|^{2}$ and $C_{2}=6k_{1}=\frac{6T^{2\alpha-1}}{\Gamma^{2}(\alpha)(2\alpha-1)}$.

Then, we consider the following integral equation:
\begin{align*}
x(t)=C_{1}+C_{2}\int_{0}^{t}L(s,x(s))ds, \tag{3.13}
\end{align*}

This equation has a globe solution via the Lemma 3.1.

Now we use the mathematical induction to prove $\|\mathbf{y}_{n}(s)\|^{2}_{L_{2}(\Omega,H)}\leq x(t), \forall~t\in [0,T]$. Firstly, we have
\begin{align*}
\|\mathbf{y}_{0}(t)\|^{2}_{L_{2}(\Omega,H)}=\mathbf{E}\|\mathbf{y}_{0}\|^{2}\leq C_{1}\leq x(t),\tag{3.14}
\end{align*}

Suppose that $\|\mathbf{y}_{n}(s)\|^{2}_{L_{2}(\Omega,H)}\leq x(t), \forall t\in [0,T]$, by using (3.12) and (3.13) and the no-decreasing property of $L$, we obtain
\begin{align*}
& x(t)-\|\mathbf{y}_{n+1}(t)\|^{2}_{L_{2}(\Omega,H)}\\
& \geq C_{2}\int_{0}^{t}(L(s,x(s))-L(s,\|\mathbf{y}_{n}(s)\|^{2}_{L_{2}(\Omega,H)}))ds\\
&\geq0.\tag{3.15}
\end{align*}

Particularly, we have $\mathrm{sup}_{n\geq0}\|\mathbf{y}_{n}\|_{L_{2}(\Omega,H)}\leq [x(T)]^{1/2}$, the Lemma 3.3 is proved.

\textbf{Lemma 3.4.} The sequence of stochastic process $\{\mathbf{y}_{n}\}_{n\geq0}$ is a Cauchy sequence.

\textbf{Proof.} Using the same argument in Lemma 3.3, we can obtain
\begin{align*}
\|\mathbf{y}_{m}(t)-\mathbf{y}_{n}(t)\|^{2}_{L_{2}(\Omega,H)}
\leq C_{3}\int_{0}^{t}K(s,\|\mathbf{y}_{m-1}(s)-\mathbf{y}_{n-1}(s)\|^{2}_{L_{2}(\Omega,H)})ds,\tag{3.16}
\end{align*}
in which $C_{3}=\frac{4T^{2\alpha-1}}{\Gamma^{2}(\alpha)(2\alpha-1)}$. Let $\rho_{n}(t)=\mathrm{sup}_{m\geq n}(\|\mathbf{y}_{m}-\mathbf{y}_{n}\|_{L_{2}(\Omega,H)}^{2})$, we imply that
\begin{align*}
\rho_{n}(t)\leq C_{3}\int_{0}^{t}K(s,\rho_{n-1}(t))ds.\tag{3.17}
\end{align*}

It is obvious that the function $\rho_{n}(t),n\geq0$ is well defined and bounded by lemma 3.3 and also monotone non-decreasing. So there exist a monotone non-decreasing function $\rho(t)$ such that $\mathrm{lim}_{n\rightarrow\infty}\rho_{n}(t)=\rho(t)$.

Using the Lebesgue convergence theorem and taking $n\rightarrow+\infty$ in above inequality, we get
\begin{align*}
\rho(t)
\leq \gamma\int_{0}^{t}K(s,\rho(t))ds.\tag{3.18}
\end{align*}

It means that $\rho(t)=0$ follows from Lemma 3.1, for all $t\in[0,T]$. However, we can see that $0\leq\|\mathbf{y}_{m}-\mathbf{y}_{n}\|_{L_{2}(\Omega,H)}^{2}\leq\rho_{n}(T)$ and $\rho_{n}(T)\rightarrow\rho(T)=0$ when $n\rightarrow+\infty$. So $\{\mathbf{y}_{n}\}_{n\geq0}$ is a Cauchy sequence.

\textbf{Theorem 3.1.} Under the conditions (3.1),(3.2),(3.4) and (3.5), there exists a unique solution of equation (2.9).

\textbf{Proof.}(1) Existence: If we denote $\mathbf{y}(t)$ by the limit of the sequence $\{\mathbf{y}_{n}\}_{n\geq0}$, repeating the proof of Lemma 3.4, then we know that the right side of second Picard's iteration (2.14) tends to
\begin{align*}
&\mathbf{y}_{0}+\frac{1}{\Gamma(\alpha)}\int_{0}^{t}(t-s)^{\alpha-1}f(s,\mathbf{y}(s))ds\\
&+\frac{1}{\Gamma(\alpha)}\int_{0}^{t}(t-s)^{\alpha-1}\sigma(s,\mathbf{y}(s))dW(s) , \tag{3.19}
\end{align*}
which is just a solution of equation (2.10).

(2) Uniqueness: Suppose $\mathbf{x}(t)$ and $\mathbf{y}(t)$ are two solutions of equation (2.10), using the same argument as in Lemma 3.3, we have
\begin{align*}
&\|\mathbf{y}(t)-\mathbf{x}(t)\|^{2}_{L_{2}(\Omega,H)}\\
&\leq C_{3}\int_{0}^{t}K(s,\|\mathbf{y}(s)-\mathbf{x}(s)\|^{2}_{L_{2}(\Omega,H)})ds, \tag{3.20}
\end{align*}

Using the Lemma 3.1 again, we can obtain $\|\mathbf{y}(t)-\mathbf{x}(t)\|^{2}_{L_{2}(\Omega,H)}=0$ for all $t\in [0,T]$, which implies that $\mathbf{y}(t)=\mathbf{x}(t)$. The proof is completed.

\section{The chaotic systems and numerical simulations}

For the interval $[0,T]$, the discretization and equidistant grid is chosen as following:
\begin{align*}
0=t_{0}<t_{1}<t_{2}<\cdots<t_{N+1}=T, t_{j+1}-t_{j}=h.\tag{4.1}
\end{align*}

The Adams predictor-corrector scheme [20,42] is modified and used to solve the stochastic fractional order differential equation (2.1), then it can be discretized as follows:
\begin{align*}
y_{h}(t_{n+1})=&\mathop{\Sigma}_{k=0}^{[\alpha]-1}\frac{t_{n+1}^{k}}{k!}\mathbf{y}^{(k)}(0)+\frac{h^{\alpha}}{\Gamma(\alpha+2)}f(t_{n+1},y_{h}^{p}(t_{n+1}))\\
&+\frac{h^{\alpha}}{\Gamma(\alpha+2)}\mathop{\Sigma}_{j=0}^{n}a_{j,n+1}f(t_{j},y_{h}(t_{j}))\\
&+\frac{h^{\alpha-1}}{\Gamma(\alpha+2)}\sigma(t_{n+1},y_{h}^{p}(t_{n+1}))\Delta W_{n}\\
&+\frac{h^{\alpha-1}}{\Gamma(\alpha+2)}\mathop{\Sigma}_{j=0}^{n}a_{j,n+1}\sigma(t_{j},y_{h}(t_{j}))\Delta W_{n},\tag{4.2}
\end{align*}
where $\Delta W_{n}=W(t_{n+1})-W(t_{n})$ denotes the Wiener increments, and
\begin{align*}
a_{j,n+1}=\begin{cases}
n^{\alpha+1}-(n-\alpha)(n+1)^{\alpha},  ~j=0, \\
(n-j+2)^{\alpha+1}-(n-j)^{\alpha+1}-2(n-j+1)^{\alpha+1}, ~1\leq j\leq n,\\
1, ~j=n+1,
\end{cases} \tag{4.3}
\end{align*}
and the predicted value $y_{h}(t_{n+1})$ is determined by the fractional Adams-Bashforth method
\begin{align*}
y_{h}^{p}(t_{n+1})=&\mathop{\Sigma}_{k=0}^{[\alpha]-1}\frac{t_{n+1}^{k}}{k!}\mathbf{y}^{(k)}(0)+\frac{1}{\Gamma(\alpha)}\mathop{\Sigma}_{j=0}^{n}b_{j,n+1}f(t_{j},y_{h}(t_{j}))\\
&+\frac{1}{\Gamma(\alpha)h}\mathop{\Sigma}_{j=0}^{n}b_{j,n+1}\sigma(t_{j},y_{h}(t_{j}))\Delta W_{n},\tag{4.4}
\end{align*}
in which $b_{j,n+1}=\frac{h^{\alpha}}{\alpha}[(n+1-j)^{\alpha}-(n-j)^{\alpha}]$.

The following stochastic fractional-order Newton-Leipnik and Lorenz system are considered, because we often make their analytical solutions impossible, so these two stochastic fractional-order chaotic systems are solved by using the above mentioned method.

\subsection{The stochastic fractional-order Newton-Leipnik system}

In 1981, Leipnik and Newton [43] found two strange attractors in rigid body motion, which is a very interesting chaotic phenomenon. Here,
the stochastic fractional-order Newton-Leipnik system is considered, and the governing equation is given as follows:
\begin{align*}
\begin{cases}
^{C}D_{t}^{\alpha}x=-\beta x+y+10yz+\sigma(x)\dot{w}_{1},\\
^{C}D_{t}^{\alpha}y=-x-0.4y+5xz+\sigma(y)\dot{w}_{2},\\
^{C}D_{t}^{\alpha}z=\rho z-5xy+\sigma(z)\dot{w}_{3},\\
\end{cases} \tag{4.5}
\end{align*}
where $\beta$ and $\rho$ are positive parameters, and usually the interval of parameter $\rho$ is taken in $[0,8.0]$. The function $\sigma(t)=\mu t$, $\mu$ is constant, and $\dot{w}_{i}=\frac{dW_{i}}{dt}, i=1,2,3$ describes a state dependent random noise. The parameters are taken as $\beta=0.4,\rho=0.175$, the initial conditions are chosen as $[0.19 \ 0 \ -0.18]$.

The differential equation systems (4.5) and the initial conditions can be written as
\begin{align*}
\begin{cases}
^{C}D_{t}^{\alpha}\mathbf{x}(t)=F(\mathbf{x}(t))+\sigma(\mathbf{x}(t))\dot{W}(t),\\
\mathbf{x}(0)=\mathbf{x}_{0},\\
\end{cases}\tag{4.6}
\end{align*}
in which
\begin{align*}
\mathbf{x}(t)=(x_{1}(t),x_{2}(t),x_{3}(t))^{T}\in R^{3},
\end{align*}
\begin{align*}
\mathbf{x}_{0}=(x_{10},x_{20},x_{30})^{T};
\end{align*}
\begin{align*}
F(\mathbf{x}(t))=A\mathbf{x}(t)+x_{2}(t)B\mathbf{x}(t)+x_{3}(t)C\mathbf{x}(t),
\end{align*}
where
\begin{align*}
A=\begin{pmatrix} -\beta & 1 & 0 \\ -1 & -0.4 & 0 \\ 0 & 0 & \rho \end{pmatrix},
B=\begin{pmatrix} 0 & 0 & 0 \\ 0 & 0 & 0 \\ 5 & 0 & 0 \end{pmatrix},
C=\begin{pmatrix} 0& 10 & 0 \\ 5 & 0 & 0 \\ 0 & 0 & 0 \end{pmatrix};
\end{align*}
\begin{align*}
\sigma(\mathbf{x}(t))=(\mu,0,0)\mathbf{x}(t)D_{1}+(0,\mu,0)\mathbf{x}(t)D_{2}+(0,0,\mu)\mathbf{x}(t)D_{3},
\end{align*}
where
\begin{align*}
D_{1}=\begin{pmatrix} 1 & 0 & 0 \\ 0 &  0& 0 \\ 0 & 0 & 0 \end{pmatrix},
D_{2}=\begin{pmatrix} 0 & 0 & 0 \\ 0 &  1 & 0 \\ 0 & 0 & 0 \end{pmatrix},
D_{3}=\begin{pmatrix} 0 & 0 & 0 \\ 0 &  0 & 0 \\ 0 & 0 & 1 \end{pmatrix};
\end{align*}
\begin{align*}
\dot{W}(t)=(\dot{w}_{1},\dot{w}_{2},\dot{w}_{3})^{T}.
\end{align*}

Obviously, $F(\mathbf{x}(t))$ and $\sigma(\mathbf{x}(t))$ are continuous and bounded on the interval $[\mathbf{x}_{0}-\delta,\mathbf{x}_{0}+\delta]$ for any $\delta>0$, furthermore, we have
\begin{align*}
&\mathbf{E}(\|F(\mathbf{x}(t))-F(\mathbf{y}(t))\|^{2})+\mathbf{E}(\|\sigma(\mathbf{x}(t))-\sigma(\mathbf{y}(t))\|^{2})\\
&=\mathbf{E}(\|A(\mathbf{x}(t)-\mathbf{y}(t))+(x_{2}(t)B\mathbf{x}(t)-y_{2}(t)B\mathbf{y}(t))\\
&\hspace{4mm}+(x_{3}(t)C\mathbf{x}(t)-y_{3}(t)C\mathbf{y}(t))\|^{2})+\mathbf{E}(\|(\mu,0,0)(\mathbf{x}(t)-\mathbf{y}(t))D_{1}\\
&\hspace{4mm}+(0,\mu,0)(\mathbf{x}(t)-\mathbf{y}(t))D_{2}+(0,0,\mu)(\mathbf{x}(t)-\mathbf{y}(t))D_{3}\|^{2})\\
&\leq \mathbf{E}((\|A\|^{2}+\|B\|^{2}(\|\mathbf{x}(t)\|^{2}+|y_{2}(t)|^{2})\\
&\hspace{4mm}+\|C\|^{2}(\|\mathbf{x}(t)\|^{2}+|y_{3}(t)|^{2}))\|\mathbf{x}(t)-\mathbf{y}(t)\|^{2})\\
&\hspace{4mm}+\mathbf{E}(\mu(\|D_{1}\|^{2}+\|D_{2}\|^{2}+\|D_{3}\|^{2})\|\mathbf{x}(t)-\mathbf{y}(t)\|^{2})\\
&\leq K_{1}\mathbf{E}(\|\mathbf{x}(t)-\mathbf{y}(t)\|^{2}), \tag{4.7}
\end{align*}
where $K_{1}=\|A\|^{2}+(\|B\|^{2}+\|C\|^{2})(2\|\mathbf{x}_{0}\|^{2}+\delta)+3\mu$.

The above inequality manifests that $F(\mathbf{x}(t))$ and $\sigma(\mathbf{x}(t))$ satisfies some non-Lipschitz conditions (A1) and (A2). Based on the results lemma 1 and lemma 2, we can
conclude that the stochastic fractional-order Newton-Leipnik system has a unique solution.

In numerical simulations, the time step is taken $h=0.005$. Fig.1 and Fig.2 give the phase portraits of stochastic fractional-order Newton-Leipnik system in $x-y-z$ space and $x-y$, $y-z$, $x-z$ planes, where the orders are taken as $\alpha=0.93$ and $\alpha=0.99$, respectively.

\setlength{\unitlength}{0.1cm}
\begin{center}
\begin{picture}(85,10)
\put(1,1){\line(1,0){80}} \put(1,1){\line(0,1){7}}
\put(1,8){\line(1,0){80}} \put(81,1){\line(0,1){7}}
\end{picture}\\
{{\bf Fig.1} }
\end{center}

\setlength{\unitlength}{0.1cm}
\begin{center}
\begin{picture}(85,10)
\put(1,1){\line(1,0){80}} \put(1,1){\line(0,1){7}}
\put(1,8){\line(1,0){80}} \put(81,1){\line(0,1){7}}
\end{picture}\\
{{\bf Fig.2} }
\end{center}

\subsection{The stochastic fractional-order Lorenz system}

The Lorenz system is a system of ordinary differential equations first studied by Edward Lorenz [44]. In this study, the stochastic fractional-order Lorenz system is considered and described by the following non-linear differential equations:

\begin{align*}
\begin{cases}
^{C}D_{t}^{\alpha}x=a(y-x)+\sigma(x)\dot{w}_{1},\\
^{C}D_{t}^{\alpha}y=cx-y-xz+\sigma(y)\dot{w}_{2},\\
^{C}D_{t}^{\alpha}z=xy-bz+\sigma(z)\dot{w}_{3},\\
\end{cases}\tag{4.8}
\end{align*}
where the parameters $a$ is the Prandtl number, $c$ is the Rayleigh number, and $b$ is the size of the region, which is approximated by the system.
The function $\sigma(t)=\mu t^{2}$, $\mu$ is constant, and $\dot{w}_{i}=\frac{dW_{i}}{dt}, i=1,2,3$ describes a state dependent random noise.

Here, the parameters are taken as $a=10,b=8/3,c=28$, the initial conditions are given as $[0.1 \ 0.1 \ 0.1]$, Eq.(4.8) represents the fractional order Lorenz
chaotic equation and the chaotic attractors of fractional order system also can be described.

The differential equation systems (4.8) and the initial conditions can be written as
\begin{align*}
\begin{cases}
^{C}D_{t}^{\alpha}\mathbf{x}(t)=F(\mathbf{x}(t))+\sigma(\mathbf{x}(t))\dot{W}(t),\\
\mathbf{x}(0)=\mathbf{x}_{0},\\
\end{cases}\tag{4.9}
\end{align*}
in which 
\begin{align*}
\mathbf{x}(t)=(x_{1}(t),x_{2}(t),x_{3}(t))^{T}\in R^{3},
\end{align*}
\begin{align*}
\mathbf{x}_{0}=(x_{10},x_{20},x_{30})^{T};
\end{align*}
\begin{align*}
F(\mathbf{x}(t))=A\mathbf{x}(t)+x_{1}(t)B\mathbf{x}(t),
\end{align*}
where 
\begin{align*}
A=\begin{pmatrix} -a & a & 0 \\ c & -1 & 0 \\ 0 & 0 & -b \end{pmatrix}, B=\begin{pmatrix} 0 & 0 & 0 \\ 0 & 0 & -1 \\ 0 & 1 & 0 \end{pmatrix};
\end{align*}
\begin{align*}
\sigma(\mathbf{x}(t))=(\mu x_{1}(t),0,0)\mathbf{x}(t)D_{1}+(0,\mu x_{2}(t),0)\mathbf{x}(t)D_{2}+(0,0,\mu x_{3}(t))\mathbf{x}(t)D_{3}, 
\end{align*}
where 
\begin{align*}
D_{1}=\begin{pmatrix} 1 & 0 & 0 \\ 0 &  0& 0 \\ 0 & 0 & 0 \end{pmatrix}, D_{2}=\begin{pmatrix} 0 & 0 & 0 \\ 0 &  1 & 0 \\ 0 & 0 & 0 \end{pmatrix}, D_{3}=\begin{pmatrix} 0 & 0 & 0 \\ 0 &  0 & 0 \\ 0 & 0 & 1 \end{pmatrix};
\end{align*}
\begin{align*}
\dot{W}(t)=(\dot{w}_{1},\dot{w}_{2},\dot{w}_{3})^{T}.
\end{align*}

Similarly, $F(\mathbf{x}(t))$ and $\sigma(\mathbf{x}(t))$ are continuous and bounded on the interval $[\mathbf{x}_{0}-\delta,\mathbf{x}_{0}+\delta]$ for any $\delta>0$, and we can obtain
\begin{align*}
&\mathbf{E}(\|F(\mathbf{x}(t))-F(\mathbf{y}(t))\|^{2})+\mathbf{E}(\|\sigma(\mathbf{x}(t))-\sigma(\mathbf{y}(t))\|^{2})\\
&=\mathbf{E}(\|A(\mathbf{x}(t)-\mathbf{y}(t))+(x_{1}(t)B\mathbf{x}(t)-y_{1}(t)B\mathbf{y}(t))\|^{2})\\
&\hspace{4mm}+\mathbf{E}(\|((\mu x_{1}(t),0,0)\mathbf{x}(t)-(\mu y_{1}(t),0,0)\mathbf{y}(t))D_{1}\\
&\hspace{4mm}+((0,\mu x_{2}(t),0)\mathbf{x}(t)-(0,\mu y_{2}(t),0)\mathbf{y}(t))D_{2}\\
&\hspace{4mm}+((0,0,\mu x_{3}(t))\mathbf{x}(t)-(0,0,\mu y_{3}(t)\mathbf{y}(t))D_{3}\|^{2})\\
&\leq \mathbf{E}((\|A\|^{2}+\|B\|^{2}(\|\mathbf{x}(t)\|^{2}+|y_{1}(t)|^{2}))\|\mathbf{x}(t)-\mathbf{y}(t)\|^{2})\\
&\hspace{4mm}+\mathbf{E}(\mu(\|D_{1}\|^{2}(\|\mathbf{x}(t)\|^{2}+|y_{1}(t)|^{2})+\|D_{2}\|^{2}(\|\mathbf{x}(t)\|^{2}+|y_{2}(t)|^{2})\\
&\hspace{4mm}+\|D_{3}\|^{2}(\|\mathbf{x}(t)\|^{2}+|y_{3}(t)|^{2}))\|\mathbf{x}(t)-\mathbf{y}(t)\|^{2})\\
&\leq K_{2}\mathbf{E}(\|\mathbf{x}(t)-\mathbf{y}(t)\|^{2}),\tag{4.10}
\end{align*}
where $K_{2}=\|A\|^{2}+(\|B\|^{2}+3\mu)(2\|\mathbf{x}_{0}\|^{2}+\delta)$.

Therefore, The above inequality also indicates that $F(\mathbf{x}(t))$ and $\sigma(\mathbf{x}(t))$ satisfies some non-Lipschitz conditions (A1) and (A2). On the basis of the results lemma 1 and lemma 2, it can be seen that the stochastic fractional-order Lorenz system has a unique solution.

We take the time step $h=0.005$, numerical results are illustrated in Fig. 3 and Fig. 4 for fractional order $\alpha=0.88$ and $\alpha=0.99$, respectively. The phase portraits of stochastic fractional-order Lorenz system in $x-y-z$ space and $x-y$, $y-z$, $x-z$ planes are shown.

\setlength{\unitlength}{0.1cm}
\begin{center}
\begin{picture}(85,10)
\put(1,1){\line(1,0){80}} \put(1,1){\line(0,1){7}}
\put(1,8){\line(1,0){80}} \put(81,1){\line(0,1){7}}
\end{picture}\\
{{\bf Fig.3} }
\end{center}

\setlength{\unitlength}{0.1cm}
\begin{center}
\begin{picture}(85,10)
\put(1,1){\line(1,0){80}} \put(1,1){\line(0,1){7}}
\put(1,8){\line(1,0){80}} \put(81,1){\line(0,1){7}}
\end{picture}\\
{{\bf Fig.4} }
\end{center}

Numerical simulation results of these two examples show that they can illustrate chaotic behaviors and stochastic attractors, and also with chaotic resonance. The performed simulations also clearly exhibit the stochastic effects in the chaotic system. With the increase of fractional orders while approaching towards standard order system, a pair of stochastic attractors are more stable.

\section{Conclusions}

In this paper, the existence and uniqueness of solution for the stochastic fractional-order differential equation systems are discussed. In particular, the nonlinear term is satisfied with the non-Lipschitz and linear growth conditions. The stochastic fractional-order Newton-Leipnik and Lorenz system are provided to show the application of our theory result. Finally, numerical simulation of these two examples are presented to illustrate the validity and feasibility of the modified Adams predictor-corrector scheme. It is worth to mention that chaos control and synchronization of this stochastic dynamic system are still interesting and significance problems, which should also be considered in the near future.

\section*{References}

[1] Hilfer R, editor. Application of fractional calculus in physics. New Jersey: World Scientific, 2001.

[2] Sabatier J, Agrawal O P, Machado J A T. Advances in fractional calculus. Dordrecht, The Netherlands: Springer, 2007.

[3] Khan N A, Ara A, Mahmood A. Approximate solution of time-fractional chemical engineering equations: a comparative study, Int. J. Chem. Reactor
Eng.,2010,8, Article A19.

[4] Hilfer R. Fractional diffusion based on Riemann-Liouville fractional derivatives, J. Phys. Chem. B, 2000, 104(16): 3914-3917.

[5] Oldham K B. Fractional differential equations in electrochemistry. Adv. Eng. Softw., 2010, 41(1): 9-12.

[6] Xu M, Tan W. Intermediate processes and critical phenomena: Theory, method and progress of fractional operators and their applications to modern mechanics. Sci. China Ser. G, 2006, 49(3): 257-272.

[7] Kulish V V, Lage J L. Application of fractional calculus to fluid mechanics. J. Fluid. Eng., 2002, 124(3): 803-806.

[8] Carpinteri A, Mainardi F. Fractals and fractional calculus in continuum mechanics. New York: Springer Wien, 1997.

[9] Magin R L. Fractional calculus in bioengineering. Redding: Begell House, 2006.

[10] Glockle W G, Nonenmacher T F. A fractional calculus approach to self-similar protein dynamics. Biophys. J., 1995, 68: 46-53

[11] Liu J G, Xu M Y. Study on a fractional model of viscoelasticity of human cranial bone. Chinese. J. Biomed. Eng., 2005, 24(1): 12-16.(in Chin)

[12] Best B J. Fractal physiology and chaos in medicine. New York: World Scientific Publishing Co Pte Ltd, 1993.

[13] Su H J, Xu M Y. Generalized viscoelastic model of otolith organs with fractional orders. Chinese. J. Biomed. Eng., 2001, 20(1): 46-52.(in Chin)

[14] \v{S}kovr\'{a}nek T, Podlubny I, Petr\'{a}\v{s} I. Modeling of the national economies in state-space: A fractional calculus approach. Econ. Model., 2012, 29(4): 1322-1327.

[15] Chen W C. Nonlinear dynamics and chaos in a fractional-order financial system. Chaos Solitons Fract., 2008, 36(5): 1305-1314.

[16] Bohannan G W. Application of fractional calculus to polarization dynamics in solid dielectric materials. Montana State University, 2000.

[17] Hilfer R. Experimental evidence for fractional time evolution in glass forming materials. Chem. Phys., 2002, 284(1): 399-408.

[18] Monje C A, Chen Y, Vinagre B M, et al. Fractional-order systems and controls: fundamentals and applications. Springer Science Business Media, 2010.

[19] Baleanu D, Machado J A T,  Luo A C. (Eds.). Fractional dynamics and control. Springer Science Business Media.2011.

[20] Agrawal S K, Srivastava M, Das S. Synchronization of fractional order chaotic systems using active control method. Chaos Solitons Fract., 2012, 45(6): 737-752.

[21] Monje C A, Chen Y Q, Vinagre B M, Xue D, Feliu V. Fractional-order systems and controls: fundamentals and applications, Springer-Verlag, London,
2010.

[22] Petr\'{a}\v{s} I. Fractional-order chaotic systems. In fractional-order nonlinear systems. Springer Berlin Heidelberg, 2011: 103-184.

[23] Li C, Liao X, Yu J. Synchronization of fractional order chaotic systems. Phys. Rev. E, 2003, 68(6): 067203.

[24] Petras I. Fractional-order nonlinear systems: modeling, analysis and simulation. Springer Science Business Media, 2011.

[25] Delbosco D, Rodino L.  Existence and uniqueness for a nonlinear fractional differential equation, J. Math. Anal. Appl., 1996, 204: 609-625.

[25] Kilbas A A, Srivastava H M, Trujillo J J.  Theory and applications of fractional differential equations, Elsevier, Amsterdam, 2006.

[26] Lakshmikantham V, Vatsala A S. Basic theory of fractional differential equations, Nonlinear Anal. TMA, 2008, 69(8): 2677-2682.

[27] Deng W. Smoothness and stability of the solutions for nonlinear fractional differential equations. Nonlinear Anal. TMA, 2010, 72(3): 1768-1777.

[28] Lin W. Global existence theory and chaos control of fractional differential equations. J. Math. Anal. Appl., 2007, 332(1): 709-726.

[29] {\O}ksendal B. Stochastic differential equations. Springer Berlin Heidelberg, 2003.

[30] Mao X. Stochastic differential equations and applications. Elsevier, 2007.

[31] Oksendal B. Stochastic differential equations: an introduction with applications. Springer Science Business Media, 2013.

[32] Schmalfu{\ss} B. The random attractor of the stochastic Lorenz system. Z. Angew. Math. Phys., 1997, 48: 951-975.

[33] Huang Z, Cao J, Jiang T. Dynamics of stochastic Lorenz family of chaotic systems with jump. J. Math. Chem., 2014, 52(2): 754-774.

[34] Anishchenko V S, Neiman A B, Safanova M A. Stochastic resonance in chaotic systems. J. Stat. Phys., 1993, 70(1-2): 183-196.

[35] Simiu E. Chaotic transitions in deterministic and stochastic dynamical systems: Applications of melnikov processes in engineering, physics, and neuroscience. Princeton University Press, 2014.

[36] Honerkamp J. Stochastic dynamical systems: concepts, numerical methods, data analysis. John Wiley Sons, 1993.

[37] Guckenheimer J. From data to dynamical systems. Nonlinearity, 2014, 27(7): R41.

[38] Sakthivel R, Revathi P, Ren Y. Existence of solutions for nonlinear fractional stochastic differential equations. Nonlinear Anal. TMA, 2013, 81: 70-86.

[39] Sakthivel R, Revathi P, Anthoni S M. Existence of pseudo almost automorphic mild solutions to stochastic fractional differential equations. Nonlinear Anal. TMA, 2012, 75(7): 3339-3347.

[40] Kerboua M, Debbouche A, Baleanu D. Approximate controllability of Sobolev type nonlocal fractional stochastic dynamic systems in Hilbert spaces. Abstr. Appl. Anal., 2013.

[41] Sousa E, Li C. A weighted finite difference method for the fractional diffusion equation based on the Riemann-Liouville derivative, Appl. Numer.
Math., 2011, 90: 22-37.

[42] Diethelm K, Ford N J, Freed A D. A predictor-corrector approach for the numerical solution of fractional differential equations. Nonlinear Dynam., 2002, 29(1-4): 3-22.

[43] Leipnik R B, Newton T A. Double strange attractors in rigid body motion. Phys. Lett. A, 1981, 86: 63-67.

[44] Lorenz E N. Deterministic nonperiodic flow. J. Atmos. Sci., 1963, 20(2): 130-141.

\bibliography{mybibfile}

\begin{figure}[htbp]
\begin{minipage}[t]{0.5\linewidth}
\includegraphics[height=6.5cm,width=6.5cm]{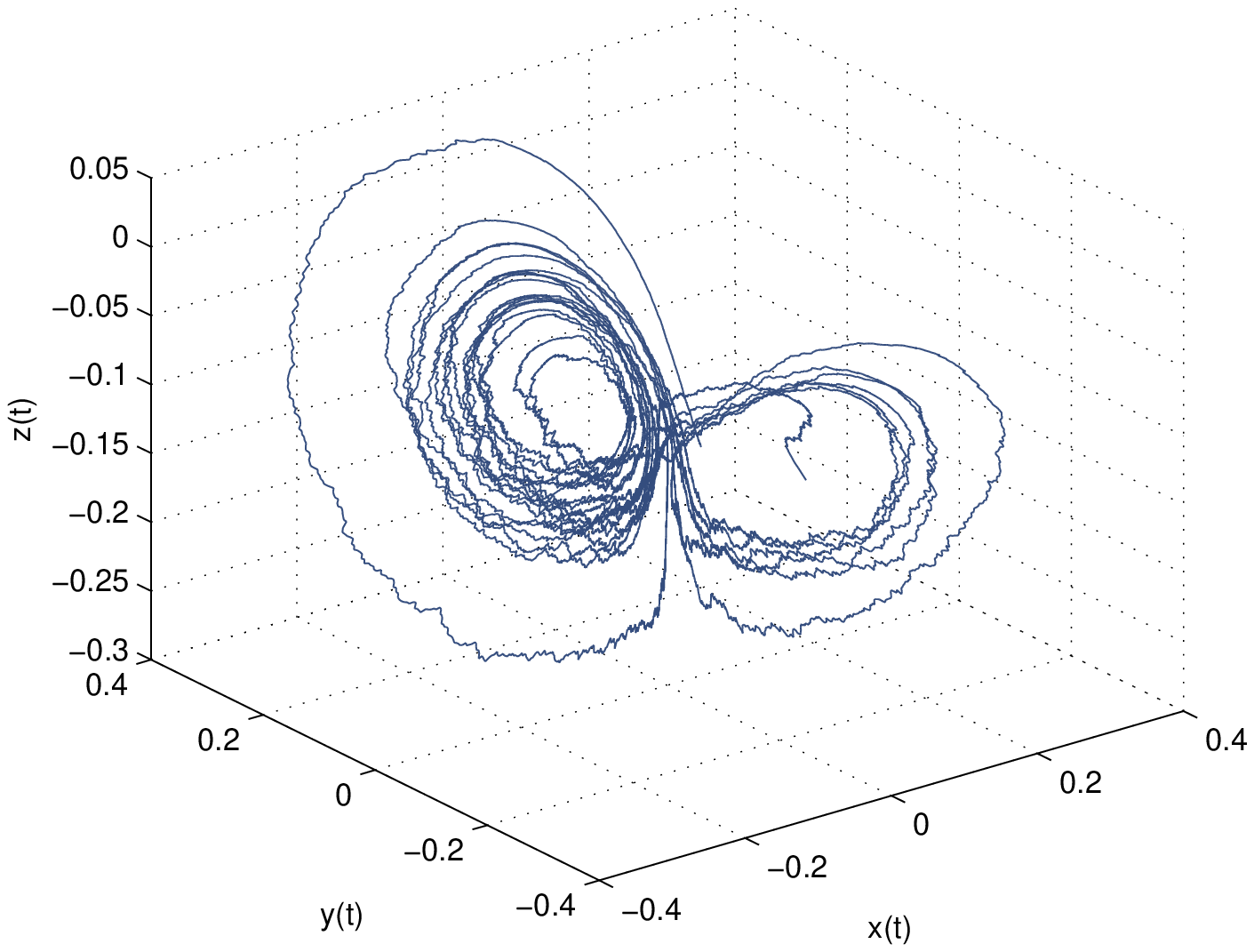}
\end{minipage}
\begin{minipage}[t]{0.4\linewidth}
\includegraphics[height=6.0cm,width=6.3cm]{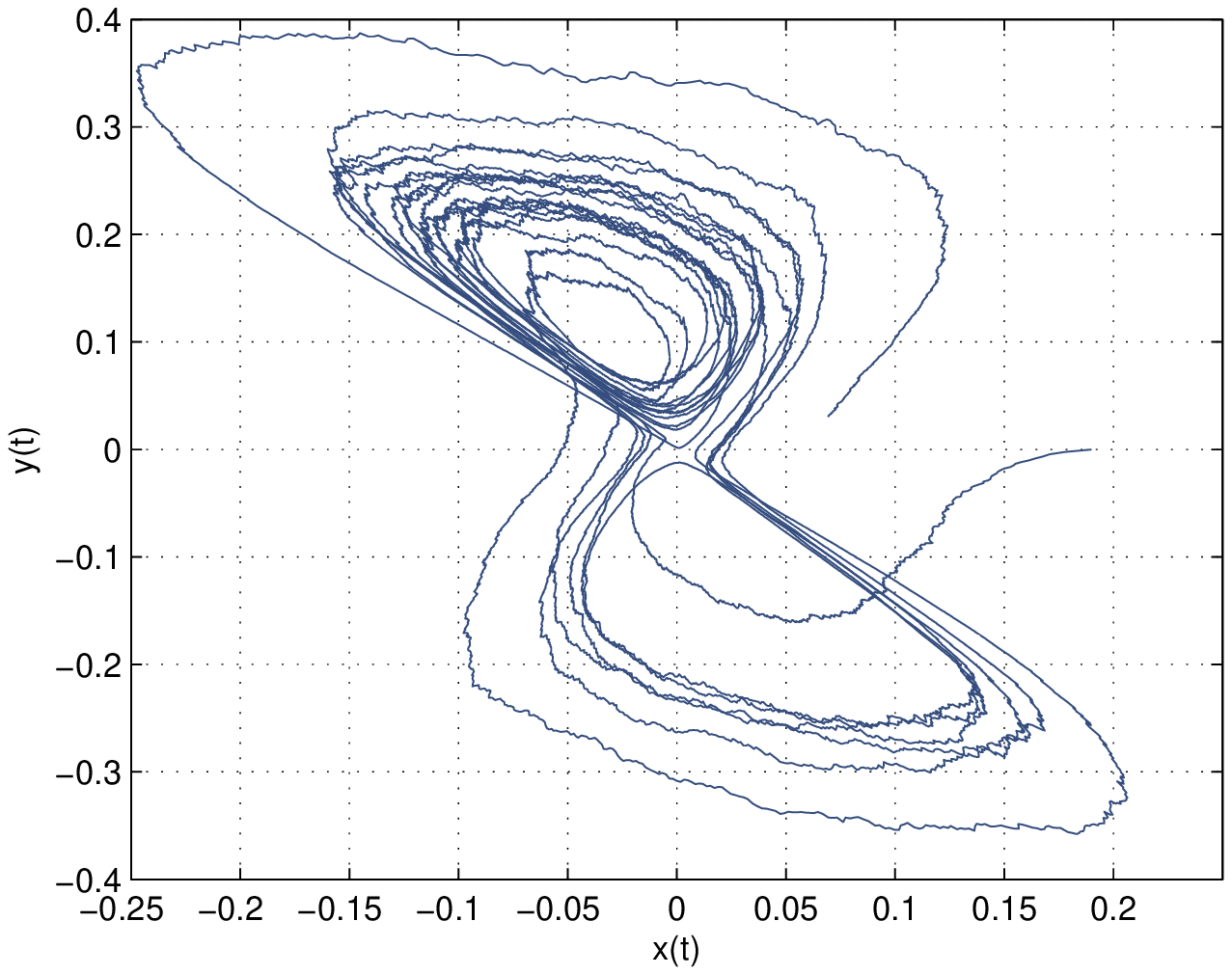}
\end{minipage}\\
\begin{minipage}[t]{0.5\linewidth}
\includegraphics[height=6.0cm,width=6.3cm]{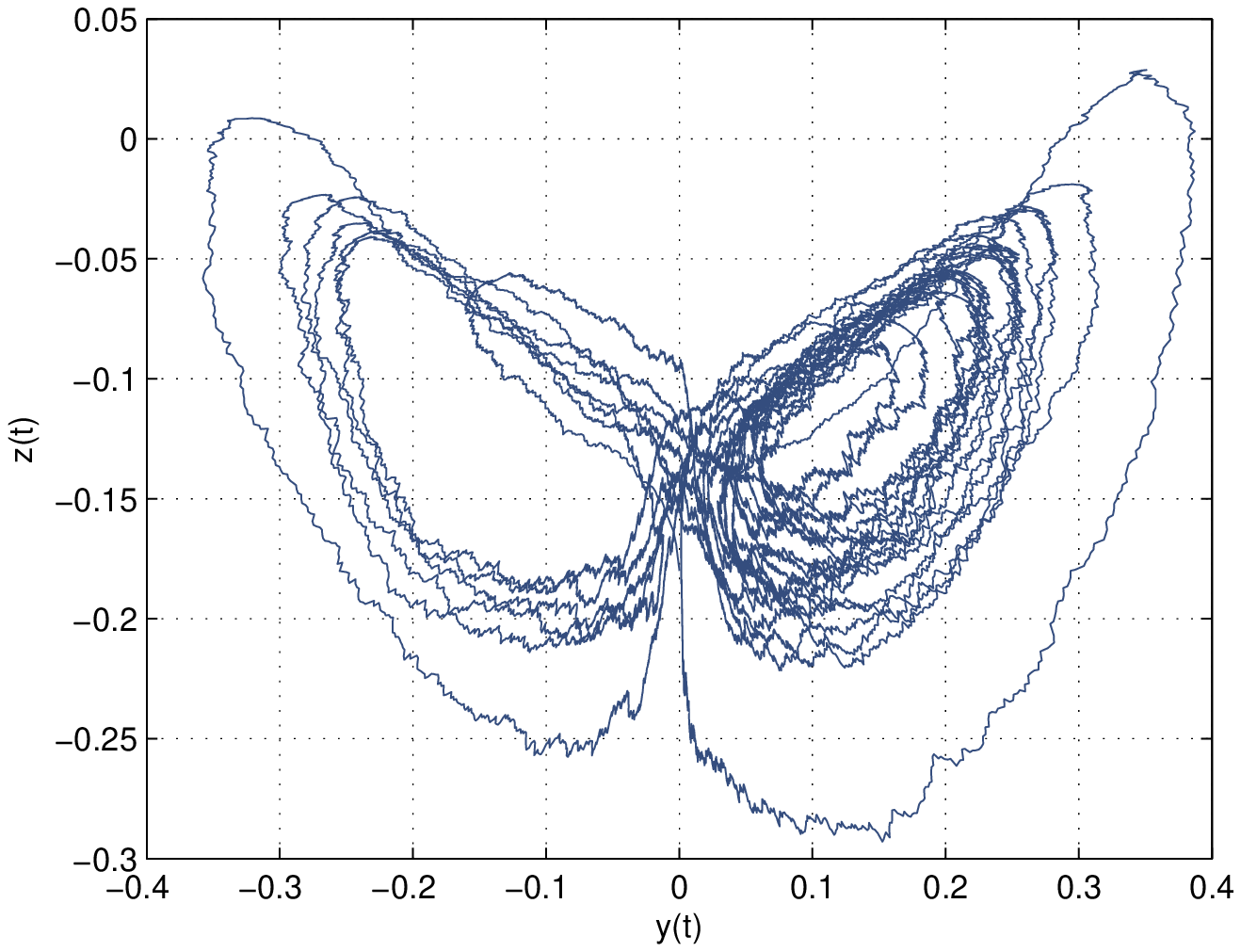}
\end{minipage}
\begin{minipage}[t]{0.4\linewidth}
\includegraphics[height=6.0cm,width=6.3cm]{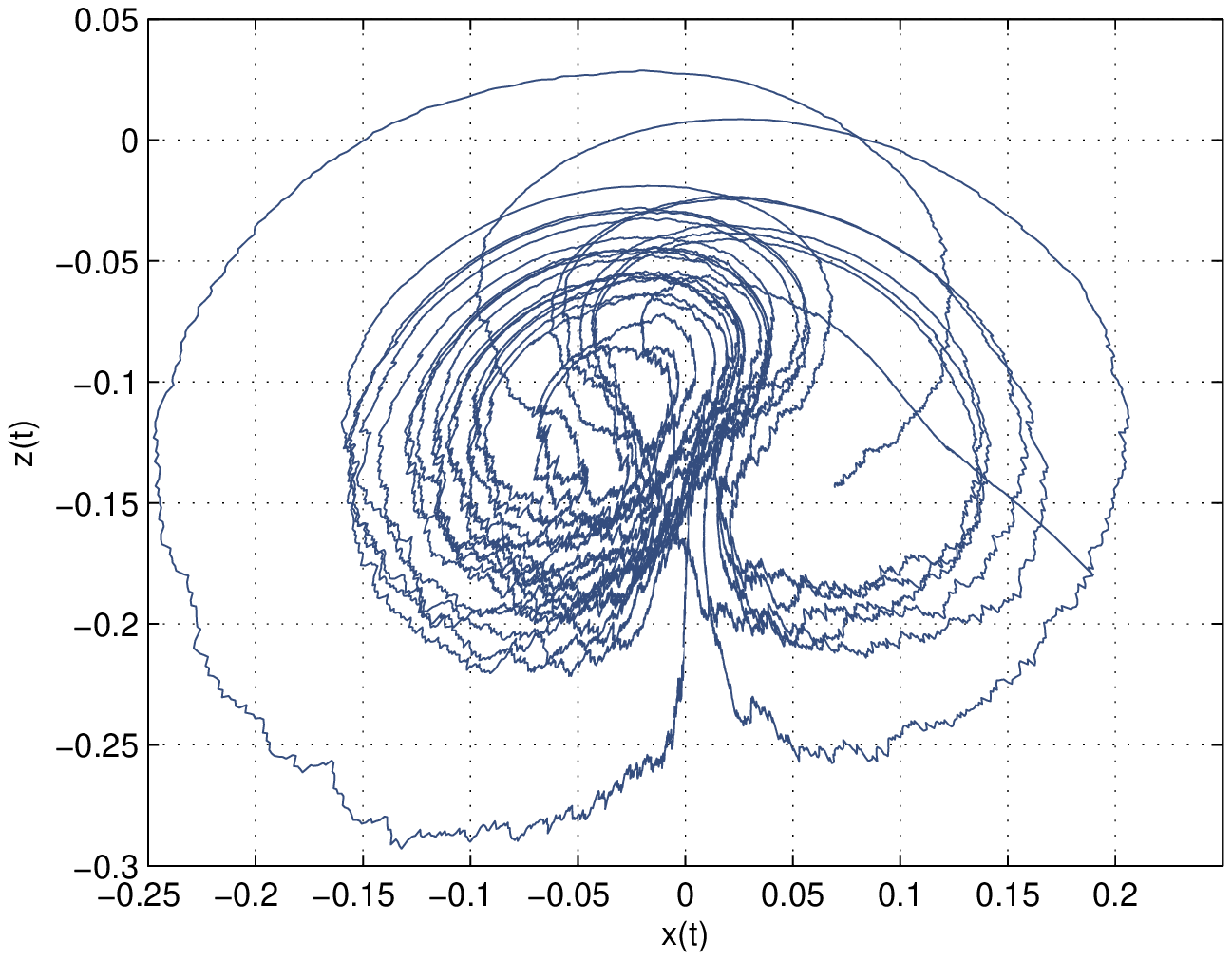}
\end{minipage}\\

{Fig.1. Phase portraits of stochastic fractional-order Newton-Leipnik system in $x-y-z$ space and $x-y$, $y-z$, $x-z$ planes. All the parameters in system are taken as $\alpha=0.93,\beta=0.4,\rho=0.175$.}
\end{figure}

\begin{figure}[htbp]
\begin{minipage}[t]{0.5\linewidth}
\includegraphics[height=6.5cm,width=6.5cm]{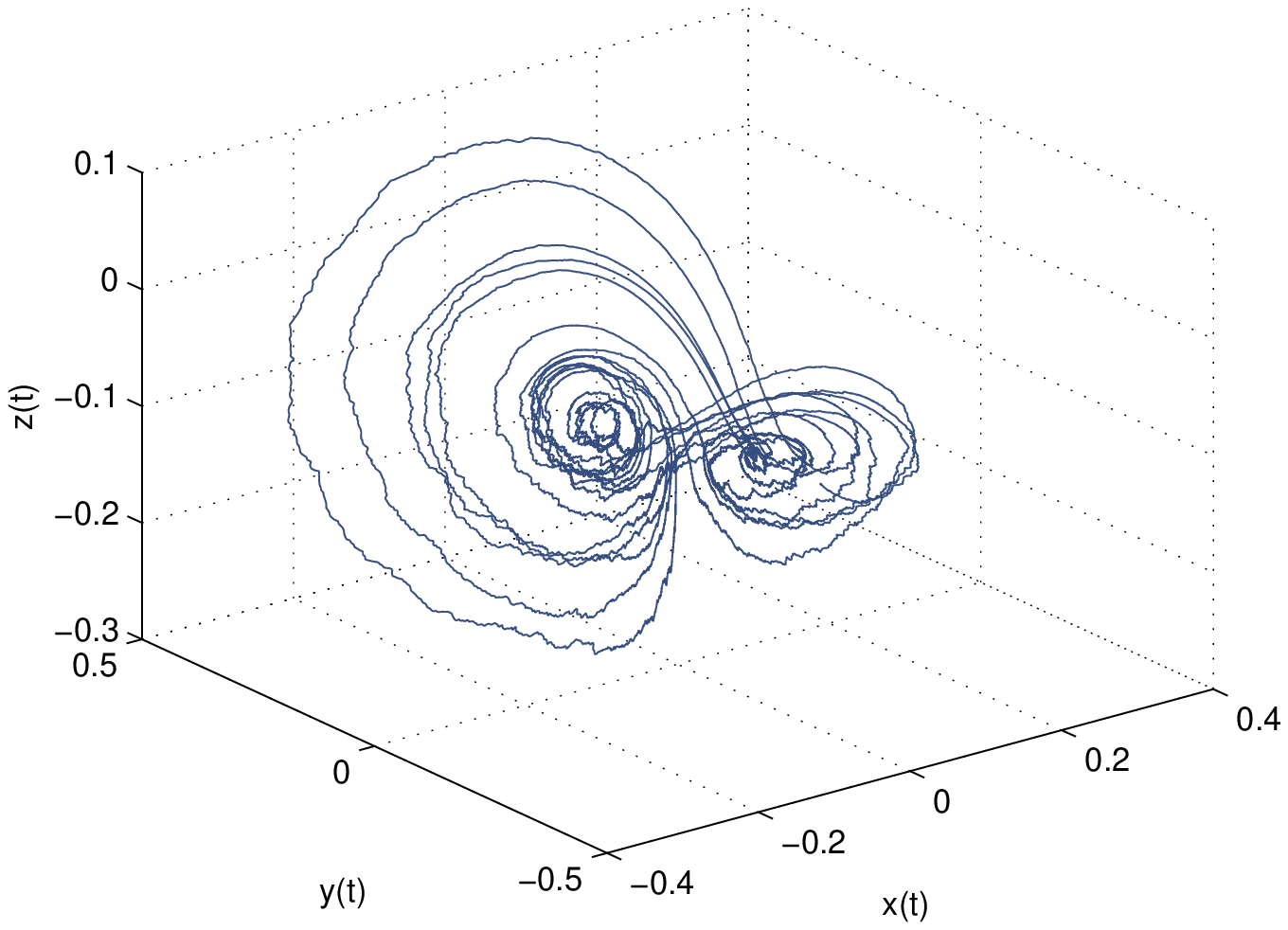}
\end{minipage}
\begin{minipage}[t]{0.4\linewidth}
\includegraphics[height=6.0cm,width=6.3cm]{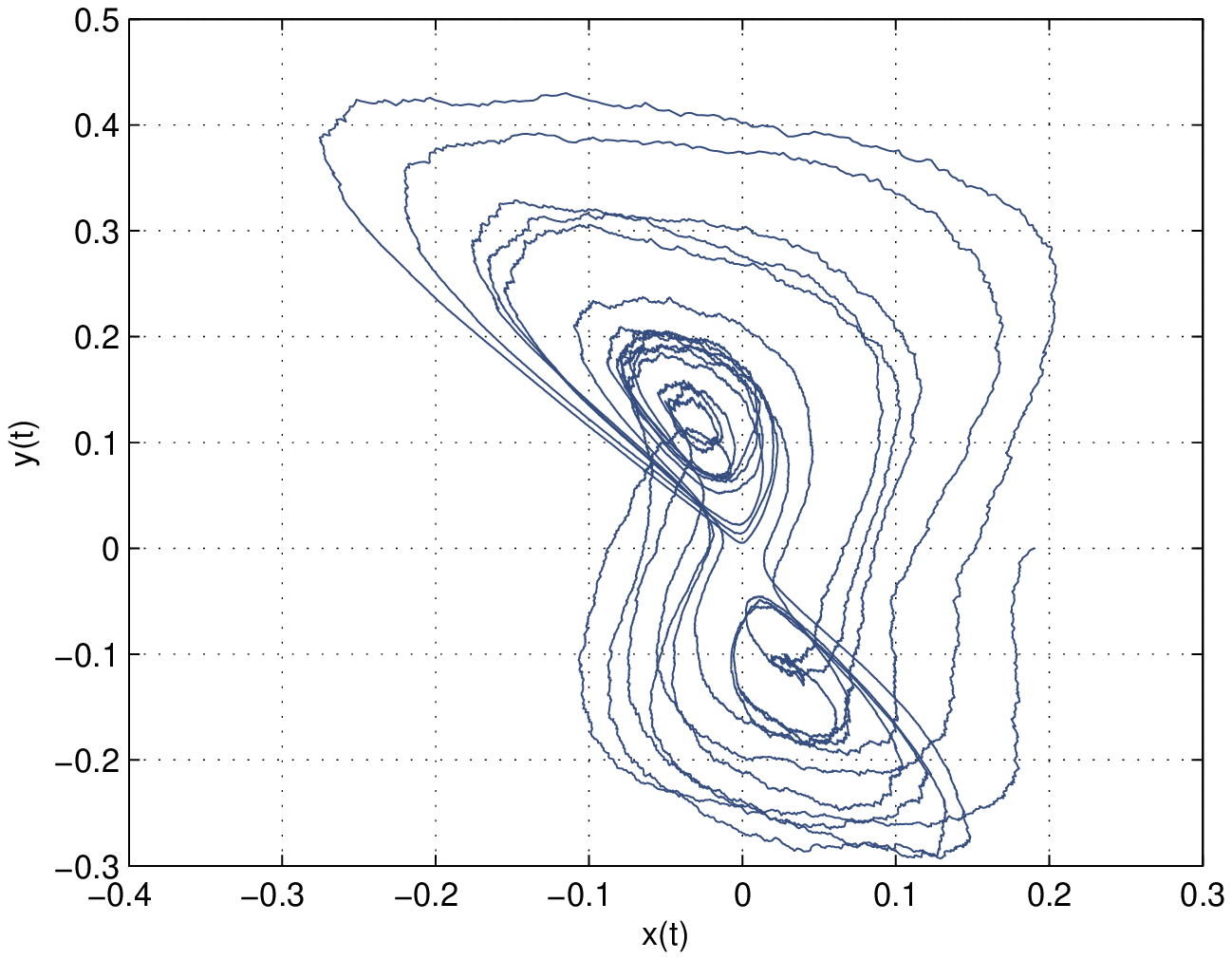}
\end{minipage}\\
\begin{minipage}[t]{0.5\linewidth}
\includegraphics[height=6.0cm,width=6.3cm]{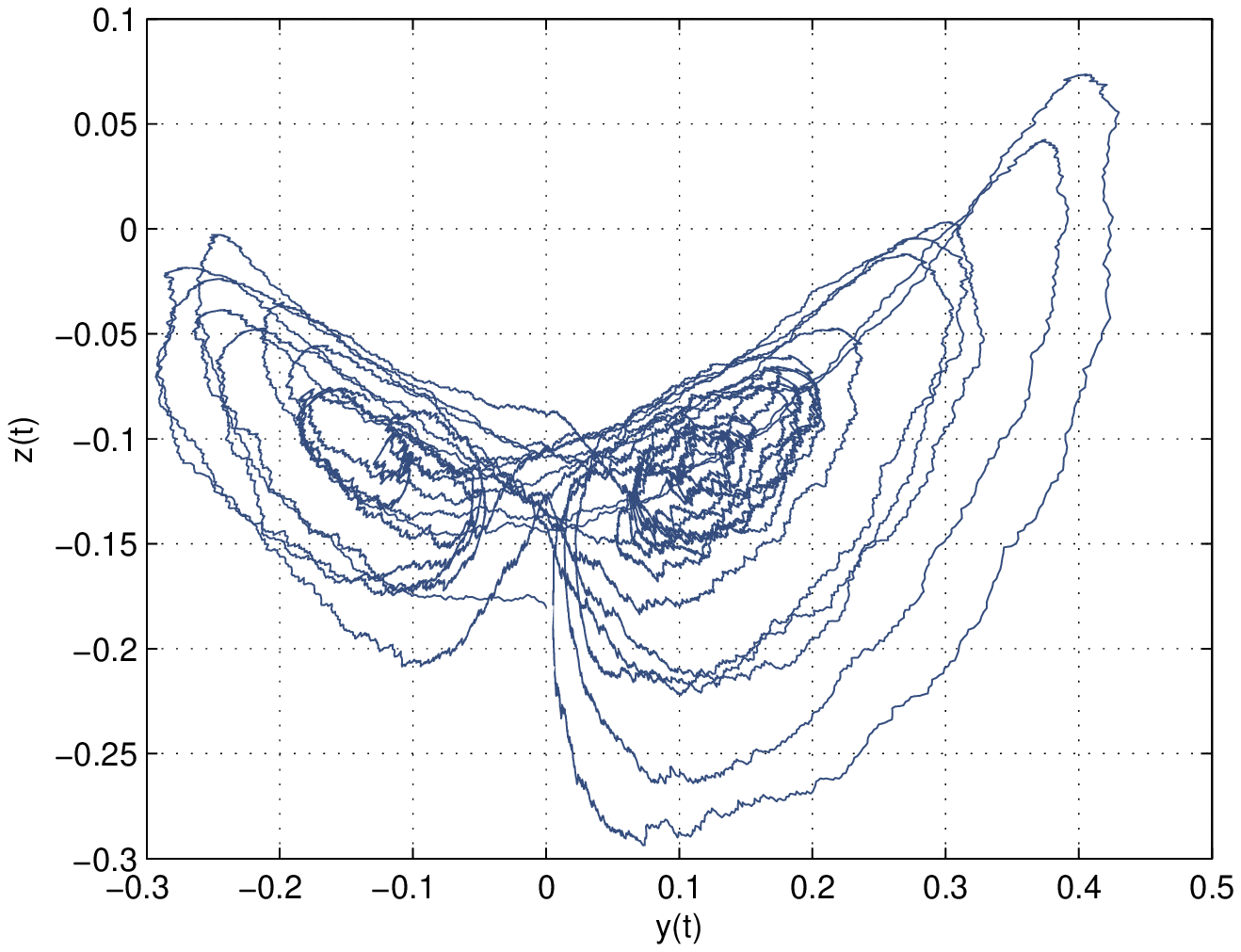}
\end{minipage}
\begin{minipage}[t]{0.4\linewidth}
\includegraphics[height=6.0cm,width=6.3cm]{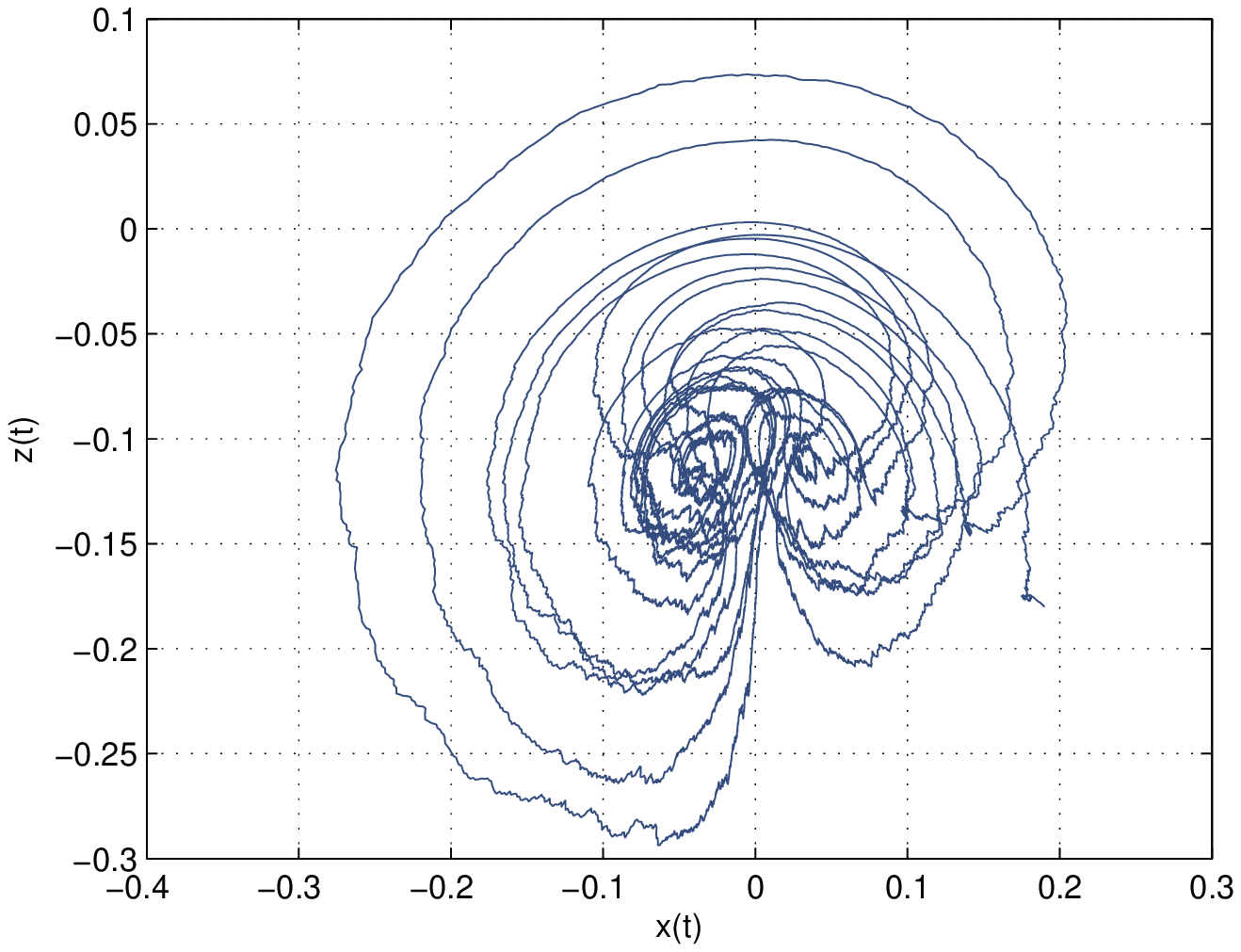}
\end{minipage}\\

{Fig.2. Phase portraits of stochastic fractional-order Newton-Leipnik system in $x-y-z$ space and $x-y$, $y-z$, $x-z$ planes. All the parameters in system are taken as $\alpha=0.99,\beta=0.4,\rho=0.175$.}
\end{figure}

\begin{figure}[htbp]
\begin{minipage}[t]{0.5\linewidth}
\includegraphics[height=6.5cm,width=6.5cm]{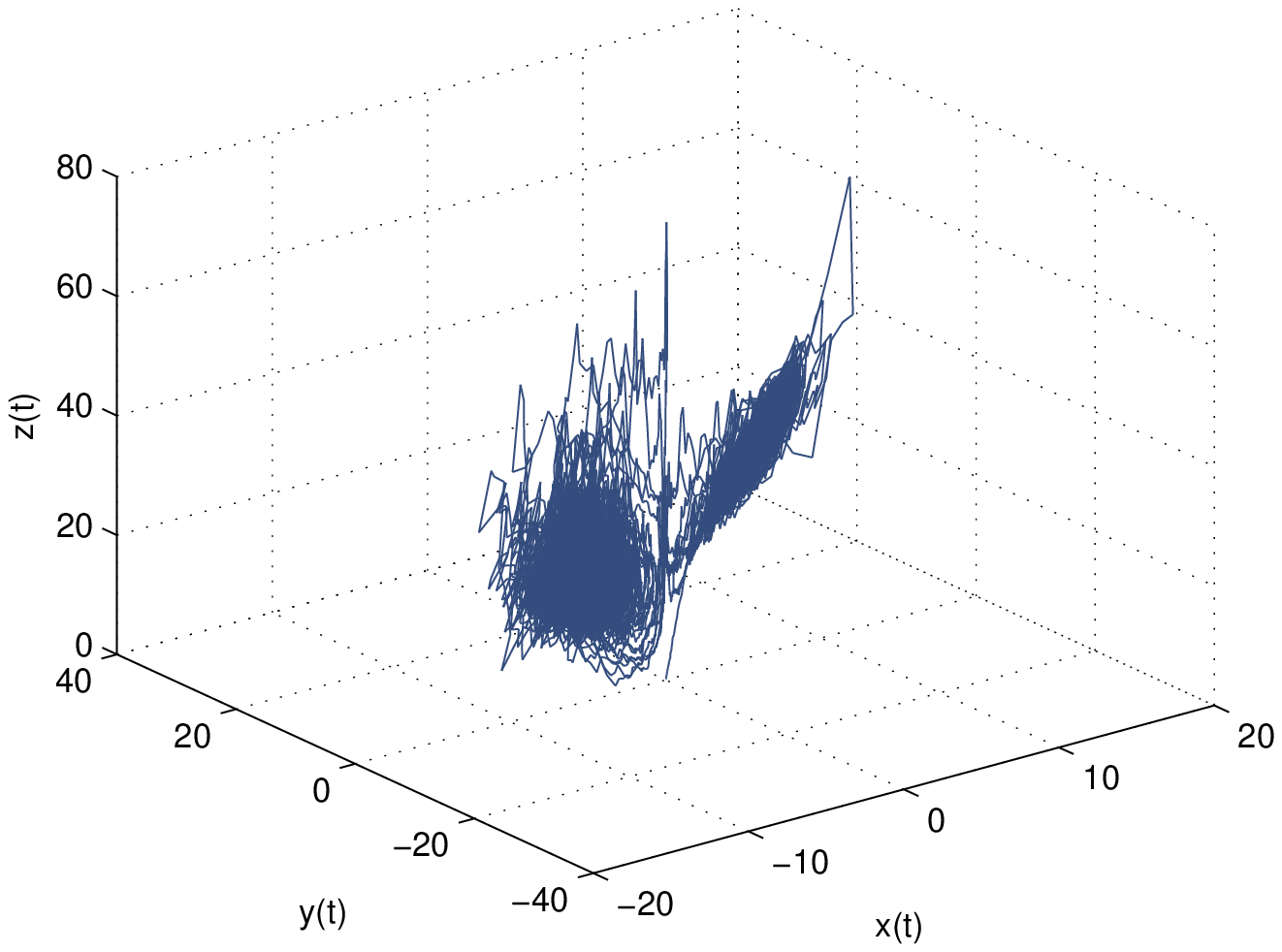}
\end{minipage}
\begin{minipage}[t]{0.4\linewidth}
\includegraphics[height=6.0cm,width=6.3cm]{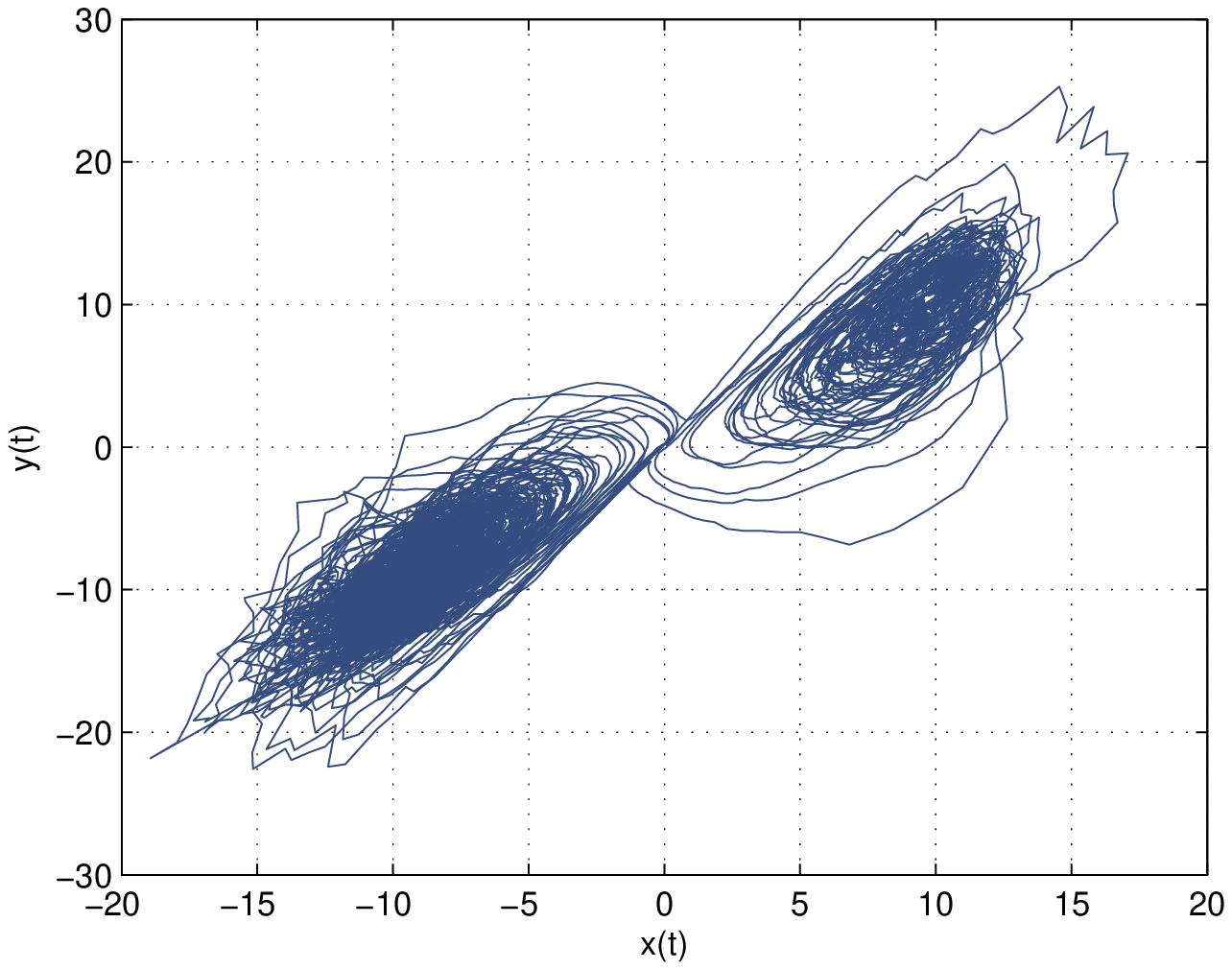}
\end{minipage}\\
\begin{minipage}[t]{0.5\linewidth}
\includegraphics[height=6.0cm,width=6.3cm]{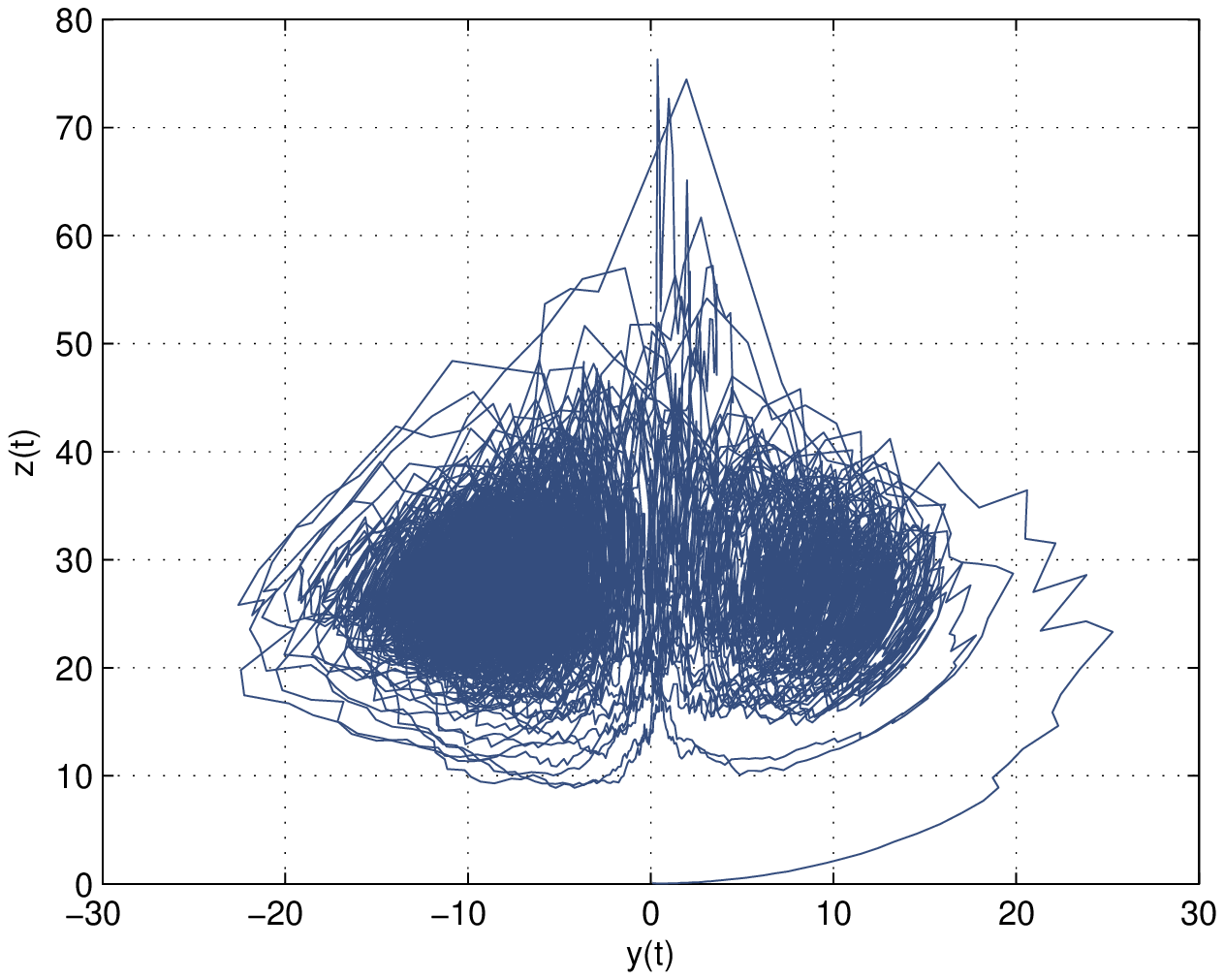}
\end{minipage}
\begin{minipage}[t]{0.4\linewidth}
\includegraphics[height=6.0cm,width=6.3cm]{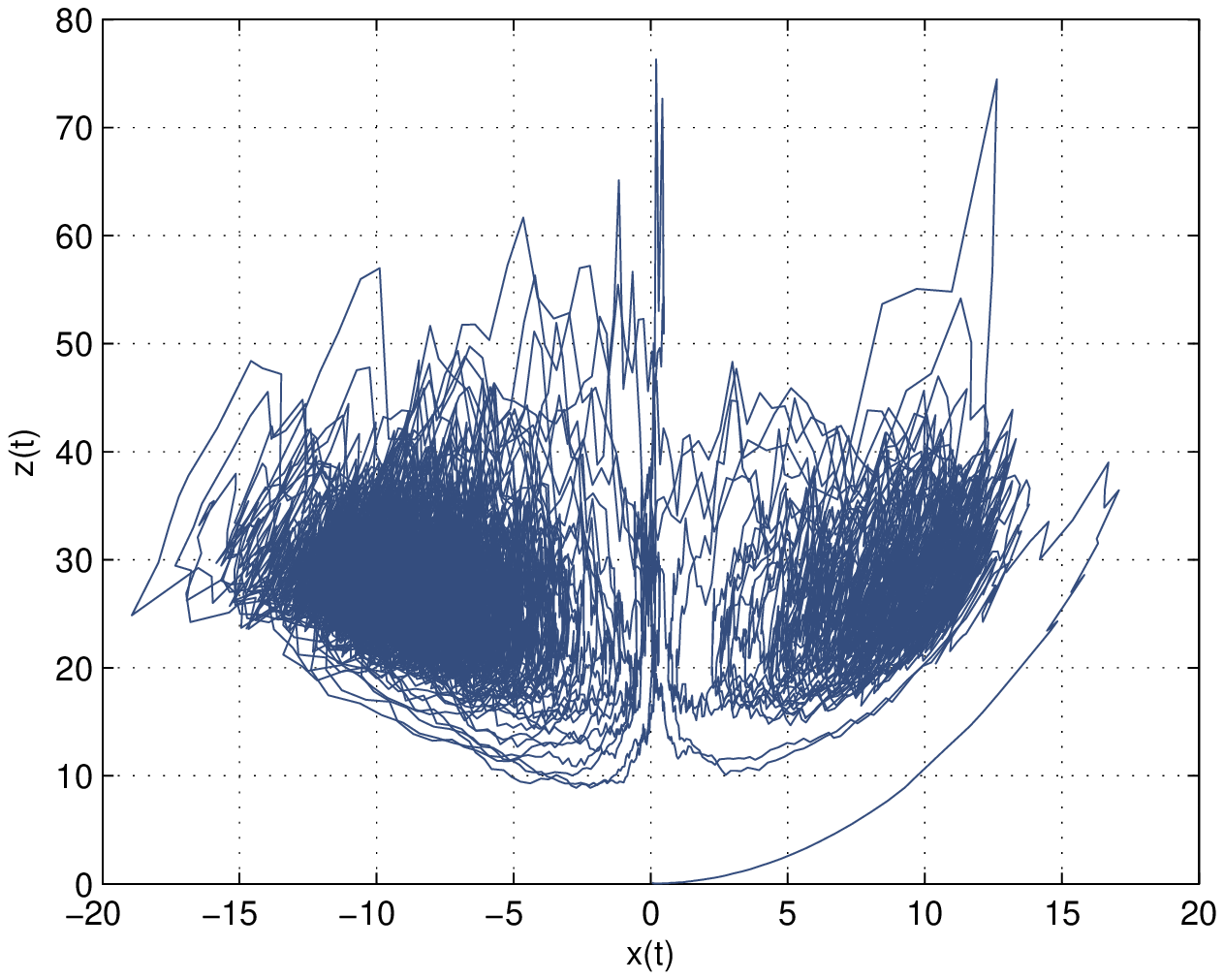}
\end{minipage}\\

{Fig.3. Phase portraits of stochastic fractional-order Lorenz system in $x-y-z$ space and $x-y$, $y-z$, $x-z$ planes. All the parameters in system are taken as $\alpha=0.88,a=10,b=8/3,c=28$.}
\end{figure}

\begin{figure}[htbp]
\begin{minipage}[t]{0.5\linewidth}
\includegraphics[height=6.5cm,width=6.5cm]{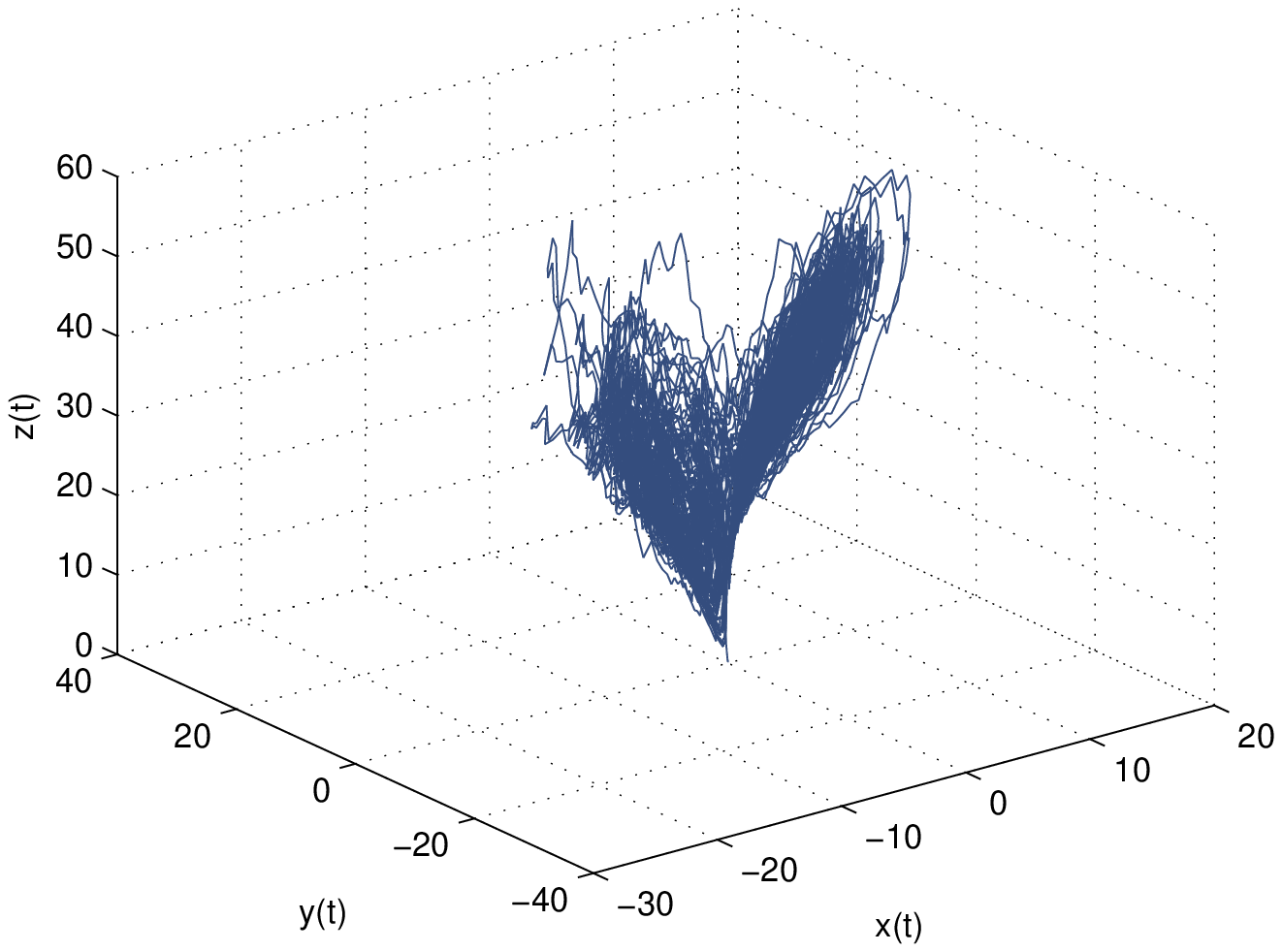}
\end{minipage}
\begin{minipage}[t]{0.4\linewidth}
\includegraphics[height=6.0cm,width=6.3cm]{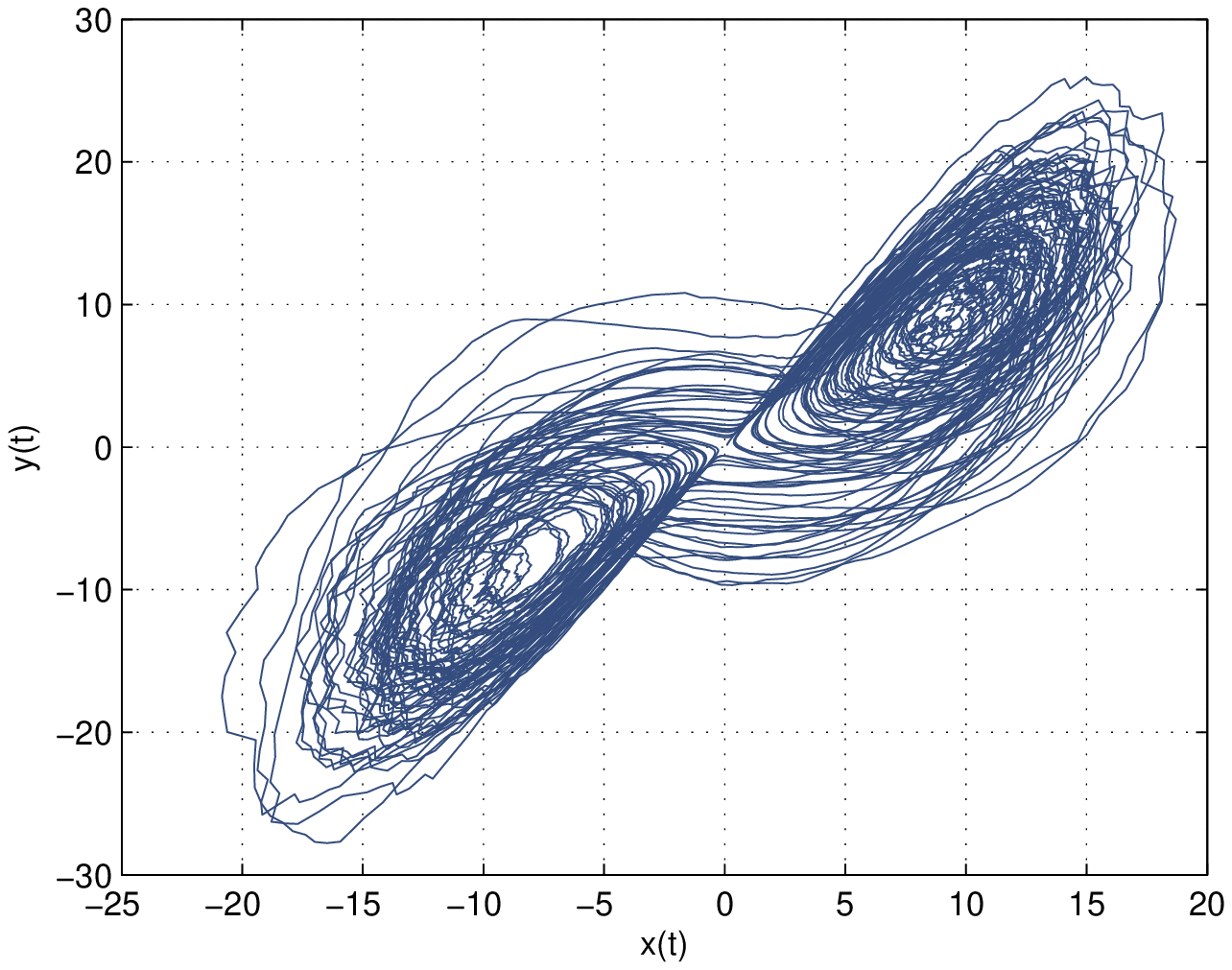}
\end{minipage}\\
\begin{minipage}[t]{0.5\linewidth}
\includegraphics[height=6.0cm,width=6.3cm]{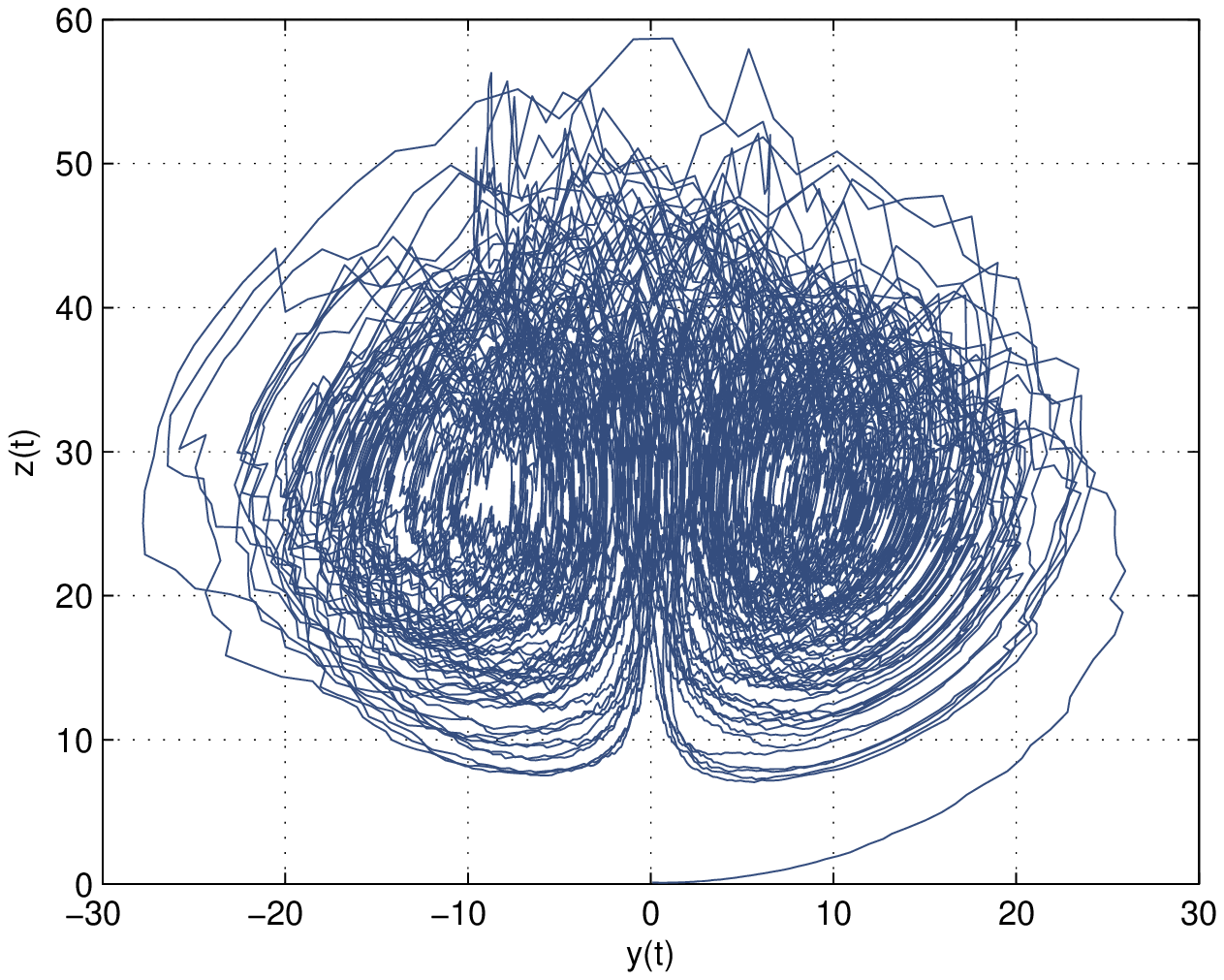}
\end{minipage}
\begin{minipage}[t]{0.4\linewidth}
\includegraphics[height=6.0cm,width=6.3cm]{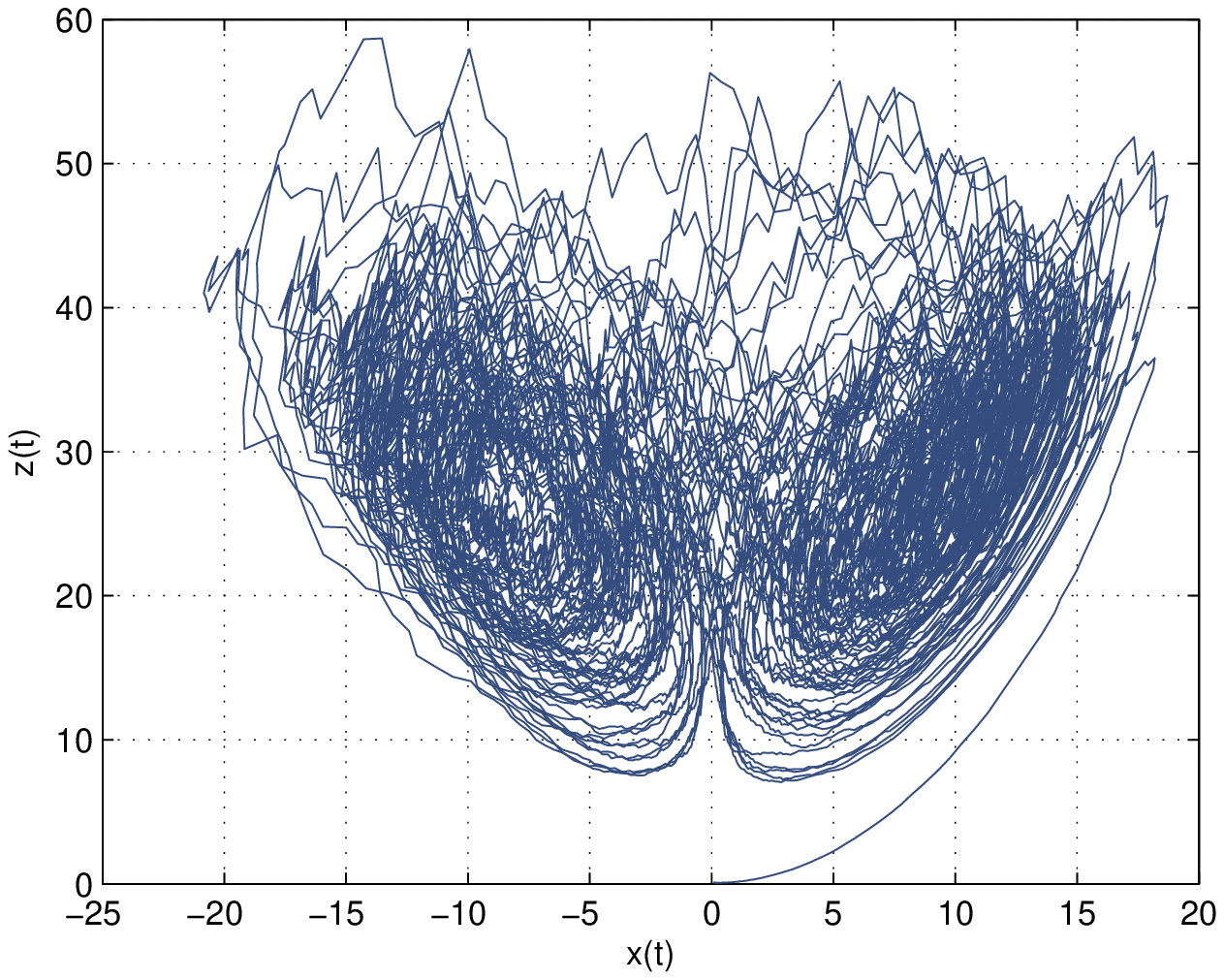}
\end{minipage}\\

{Fig.4. Phase portraits of stochastic fractional-order Lorenz system in $x-y-z$ space and $x-y$, $y-z$, $x-z$ planes. All the parameters in system are taken as $\alpha=0.99,a=10,b=8/3,c=28$.}
\end{figure}

\end{document}